\input amstex
\input psfig.sty
\documentstyle{amsppt}


\topmatter
\title
A Diameter Bound for Closed Hyperbolic 3-Manifolds.
\endtitle
\author
Matthew E. White
\endauthor
\date
Version 3.1  10 May 2000
\enddate
\endtopmatter
\document
\head 1. Introduction
\endhead
\medskip

Let $P=<a_1\dots a_n\ |\ r_1\dots r_m>$ be a presentaion for a group
$G$, then the {\bf length} of $P$ is given by $$\ell(P) = \sum_{i=1}^m
\ell(r_i)$$ where $\ell(r_i)$ is the word length of the relation $r_i$.
In [C], Cooper proves the following:

\proclaim{Theorem 1.0 (Cooper)} For every closed hyperbolic 3-manifold
$M$, the following relationship holds: $$vol(M)\ < \pi \ell(P)$$ for
every presentation $P$ of $\pi_1 (M).$ \endproclaim

In this chapter, we are going to prove a similar theorem
for the diameter of a closed, orientable hyperbolic 3-manifold:  

\proclaim {Theorem 5.9} There is an explicit constant $0 < R$ 
such that if $M$ is a closed, connected, hyperbolic 3-manifold, and $P$ is 
any presentation of its fundamental group, then $diam(M)< R(\ell(P))$.
\endproclaim

By Margulis' result, in our setting diameter and injectivty radius are 
inversely related.  Thus, our theorem can also be viewed as a lower bound 
on injectivity radius; that is, with the above hypothesis, 
$inj(M) > {1\over R(\ell(P))}$.  It is known that that infinitely many closed, 
hyperbolic 3-manifolds of volume less than a given upper bound may be 
obtained by hyperbolic Dehn surgery on a finite list of compact manifolds. 
But only finitely many of these closed manifolds have diameter less than 
a given upper bound.  Thus, our results provide a sharper version
of Theorem 1.0.

An outline of the proof is the following: we 
construct a straight 2-complex from a fundamental group presentation 
which maps into the manifold $\pi_1$-isomorphically.  
Assuming the diameter of the mainifold is very large compared to
the presentation length, Margulis provides us with a deep solid torus
surrounding a short geodesic core.  It turns out that the 2-complex cannot be homotoped
to be disjoint from the solid torus.  We consider the subcomplex that maps into the solid torus.
From this subcomplex, we construct another 2-complex which simultaneously has a large torsion 
subgroup in first homology and a small triangulation.  This is a contradiction.

Our discussion is organized as follows: In Section 2, we define the 
presentation complex and quote some of Cooper's results from his proof of the
bound for volume.  Section 3 discusses Margulis' Lemma and describes its 
influence on the geometery of the presentation complex.  Section 4 shows 
that we can construct a 2-complex with first homology bounded by 
presentation length which maps into the Margulis solid torus.  The main 
theorem is proved in Section 5.   

\head 2. The Triangular Presentation and its Complex
\endhead
\medskip

Given a presentation $P$ of a group $G$, a {\bf presentation complex} $K$ for $P$ is
a 2-complex with fundamental group $G$ constructed as follows.  The 1-skeleton $
K^{(1)}$ is a wedge of circles, one for each generator of $P$.  There is one 2-cell ${\Bbb D}$
for each relator $r$ of $P$, where $\partial {\Bbb D}$ is glued to $K^{(1)}$ along a loop
representing $r$.  If $M$ is a closed hyperbolic 3-manifold, $P$ is a presentation of 
$\pi_1(M)$, and $K$ is a presentation complex for $P$, then there is a map 
$f:K\longrightarrow M$ that induces an isomorphism of fundamental groups.

A presentation $P$ for a group $G$ is called {\bf minimal} if it has minimal length
over all presenations of $G$.  If $M$ is a closed hyperbolic 3-manifold, then $\pi_1(M)$ is
finitely generated so a minimal presentation for $\pi_1(M)$ exists. 

A presentation is called {\bf triangular} if every relation has length 3.  
We have the following:

\proclaim {Proposition 2.0} If $P$ is a finite presentation of a group $G$, then there
is a triangular presentation $P^\prime$ of $G$ such that $\ell(P^\prime)\le 3\ell(P)$.
\endproclaim

\noindent {\bf Sketch Proof} Let $P=<a_1\dots a_n\ |\ r_1\dots r_m>$. 
Proceed by induction on $\ell(P)$.  Remove all relators of the form $a_i=a_j$
by replacing every occurance of $a_j$ with $a_i$.  For a relator of the form
$a_i^2=1$, add a generator $b$ and the relations $a_i^2 b=1$, $b_i b_i b^{-1}
=1$. If there is a relator of length at least 4 $a^{e_1}_{i_1}
a^{e_2}_{i_2}\dots a^{e_k}_{i_k}$, add a generator $b$ and two relations
$a^{e_1}_{i_1} a^{e_2}_{i_2} b = 1$, 
$b^{-1} a^{e_3}_{i_3}\dots a^{e_k}_{i_k} = 1$. $\square$ 

Now, let $P=<a_1\dots a_n\ |\ r_1\dots r_m>$ be a triangular presentation
for $\pi_1(M)$. We may choose $f:K^{(1)}\longrightarrow M$ so that the image
has straight edges (each edge is a geodesic lasso at $f(v)$).
The image of each relator $r_i$ lifts to a loop 
in ${\Bbb H^3}$ with three geodesic edges.  These edges bound a hyperbolic triangle.  
Extend $f$ to map the 2-cell corresponding to $r_i$ to this triangle.  Repeat this procedure for
each relator.  It follows that $f:K\longrightarrow M$ is a $\pi_1$-isomorphism.  We call 
$f:K\longrightarrow M$ a {\bf triangular complex}. 
Clearly, each relation in $P$ corresponds to a disc in $K$ which lifts to a genuine
hyperbolic triangle in ${\Bbb H}^3$. This constuction is used by Cooper [C] to prove 
Theorem 1.0. The proof of this theorem also provides these two important
propositions: 
\proclaim {Proposition 2.1 (Bounded Area)} Using the pull-back metric,
$Area(K)\le\pi\ell (P).$ \endproclaim
\medskip

\proclaim {Proposition 2.2 (Invariant Intersection)} Suppose $K$ is a
2-complex and $f:K\longrightarrow M$ induces a $\pi_1$-isomorphism. 
Then every essential loop in $M$ meets $f(K)$.
\endproclaim
\medskip

\head 3. Geometric Preliminaries
\endhead
\medskip

In this section, we use the results above together with Margulis'
well-known and fundamental result to develop a picture of the
``thick-thin'' decompostion for a hyperbolic 3-manifold which is
relevent to our setting. A nice (and detailed) overview of this material
may be found in [BP]. We are concerned with closed, connected, oriented,
hyperbolic 3-manifolds. Therefore, unless specifically noted otherwise,
every manifold we consider below shall be of this type. By
the {\bf $\epsilon$-thick part} of $M$, we mean $$M_{[\epsilon,\infty)}
= \lbrace \text{ } x\in M \text{ }|\text{ } inj(x) \ge\epsilon\text { }
\rbrace.$$ Likewise, the {\bf
$\epsilon$-thin part} is $$M_{(0,\epsilon]} = \text{closure}
(M-M_{[\epsilon,\infty)}).$$ The structure of the thin part of a closed
hyperbolic 3-manifold is displayed in the following result due to
Margulis [BP]: 

\proclaim{Theorem 3.0 (Margulis)} Suppose that $M$ is a
closed, orientable hyperbolic 3-manifold. There exists a universal
constant $\tilde\epsilon$ such that if $\epsilon\le\tilde\epsilon$, then
$M_{(0,\epsilon]}$ is a disjoint union of solid tori (``Margulis
Tubes''). The degenerate case where one of the solid tori is 
$S^1$ may occur. Moreover, the core curve of each solid torus is a 
geodesic in $M$ which generates an infinte cyclic subgroup of $\pi_1(M)$.
\endproclaim
\medskip

We can now use Theorem 1.0 to obtain the analog of the main theorem for
the thick part of a hyperbolic manifold: 
\proclaim{Lemma 3.1 (Thick Diameter is Bounded)} There is as universal
constant $C_1$ such that given a closed hyperbolic 3-manifold  $M$ with 
$diam(M_{[\tilde\epsilon,\infty)}) \le C_1
vol(M)$ for every presentaion $P$ of $\pi_1(M)$. Thus, 
$diam(M_{[\tilde\epsilon,\infty)}) \le C_1 \ell(P)$ for every presentaion $P$ of $\pi_1(M)$ 
\endproclaim

\noindent{\bf Proof} Since $M_{[\tilde\epsilon,\infty)}$ is compact, there
exist points $x$ and $y$ in $M_{[\tilde\epsilon,\infty)}$ such that $d(x,y) =
diam(M_{[\tilde\epsilon,\infty)})$.  Now, $M_{[\tilde\epsilon,\infty)}$
is path connected, so there is a path $L$ in $M_{[\tilde\epsilon,\infty)}$
with endpoints $x$ and $y$. By considering 
$M_{[{\tilde\epsilon\over 2},\infty)}$ if necessary, we can arrange that 
each point on $L$ is conatined in a ball of radius $\tilde\epsilon\over 4$
that is isometric to a standard ball $B({\tilde\epsilon\over 4})$ in 
${\Bbb H}^3$

There is a covering of $L$ so that each ball $B$ meets at 
most two others in the covering and every such intersection is a point of
tangency between the boundary of $B$ and the boundary of another ball 
in the covering.  It follows that since the interiors of the balls do not intersect, the 
volume of the covering is the sum of the volumes of the balls in the covering.  Our
construction implies that at most ${vol(M)\over B({\tilde\epsilon\over 4})} + 2$ 
balls cover $L$. Then using Theorem 1.0, $$diam(M_{[\epsilon,\infty)})\le
2\tilde\epsilon ({vol(M)\over B({\tilde\epsilon\over 4})} + 2) \le 
{2\tilde\epsilon (\pi + 2)\over B({\tilde\epsilon\over 4})} \ell (P)$$ 
for every presentaion $P$ of $\pi_1(M)$.
Therefore, put $C_1 = {2\tilde\epsilon (\pi + 2) \over 
B({\tilde\epsilon\over 4})}$ and the proof is complete.$\square$ 
\medskip

Notice that Proposition 2.2 (invariant intersection) implies that 
if $f:K\longrightarrow M$ induces an isomorphism of fundamental groups, 
then $f(K)$ meets the core curve of every Margulis
Tube in $M$. Our approach shall be to study the intersection of $f(K)$
with a particular tube, the idea being that the geometry of a solid
torus is simple enough to impose significant contraints on the map $f$.
One can combine Theorem 1.0 and Lemma 3.1 in the following way: given
$M$, if $diam(M)$ is very large compared to $\ell(\pi_1(M))$, then
$diam(M)$ is very large compared to $vol(M)$. Both
$diam(M_{[\tilde\epsilon,\infty)})$ and $vol(M)$ are bounded above by a
constant multiple of $\ell(\pi_1(M))$. Hence, one of the Margulis Tube
components of $M_{(0,\tilde\epsilon]}$ must be have large diameter. This 
is the content of the following: 

\proclaim{Proposition 3.2 (Deep Tube)} Let $M$ be as above, let 
$C>C_1 + 1$ and suppose that $diam(M)\ge C\ell(P)$ for some
presentation $P$ of $\pi_1(M)$. Then there exists a Margulis Tube
$V_{\tilde\epsilon}\subset M$ such that $diam(V_{\tilde\epsilon}) 
\ge {1\over 2}(diam(M)- C_1 vol(M))$. Thus,
$diam(V_{\tilde\epsilon}) 
\ge {1\over 2}(C-1)\ell(P)$. \endproclaim

\noindent {\bf Proof} Given $M$, suppose that $diam(M)\ge C\ell(P)$ for 
some presentation $P$ of $\pi_1(M)$. By comapctness, there exist $x$ and 
$y$ in $M$ so that $d(x,y) \ge C\ell(P)$. 
Lemma 3.1 (Thick Diameter is Bounded) shows
that at least one of the points is contained in
$M_{(0,\tilde\epsilon]}$. Thus, if one of the points, say $y$, is
contained in the thick part of $M$, and $x\in V_{\tilde\epsilon}\subset
M_{(0,\tilde\epsilon]}$, then $d(y,\partial V_{\tilde\epsilon})\le
C_1\ell (P)$. This implies that $$(C - C_1)\ell(P)\le d(x,y)-d(y,\partial
V_{\tilde\epsilon})\le d(x,V_{\tilde\epsilon})$$ whence
$diam(V_{\tilde\epsilon}) \ge (C-1)\ell(P).$ As similar
approach shows that if both $x$ and $y$ are contained in
$M_{(0,\tilde\epsilon]}$, then: 
$$diam(V_{\tilde\epsilon}) \ge {1\over 2} (C-1)\ell(P).$$ 
$\square$ \medskip

We shall use the notation $V_{\tilde\epsilon}$ for a deep
tube of $M$ provided by the above proposition. It is useful for us to
have a picture of the geometry of $\partial V_{\tilde\epsilon}$. To
develop this, consider the upper-half space model of ${\Bbb H}^3$. Let
$\gamma$ denote the core curve of $V_{\tilde\epsilon}$. In $\pi_1(M)$,
$[\gamma]$ corresponds to a loxodromic isometry $\tilde\gamma$ of ${\Bbb
H}^3$. We may assume that the axis of $\tilde\gamma$ is the z-axis, so
that a component of the preimage of $V_{\tilde\epsilon}$ under the
universal cover $\rho:{\Bbb H}^3\longrightarrow M$, denoted $\tilde
V_{\tilde\epsilon}$, is a neighborhood of the z-axis. This appears as an
infinite Euclidean cone with (ideal) vertex at the origin. The isometry
$\tilde\gamma$ acts by translation along and rotation around the z-axis,
so that evidently a fundamental domain for $V_{\tilde\epsilon}$ is a
horizontal ``pancake;'' that is, in cylindrical coordinates the fundamental domain 
is a set of the form $\lbrace (r,\theta,z)\text{ | }a\le z\le b, {r\over z} \le c\rbrace$. 
The covering translation glues this set top to bottom with a twist. If $C>C_1$ and $diam(M)\ge
C\ell(P)$, Margulis and Proposition 3.2 (deep tube) tell us that the
radius of this ``pancake'' is at least ${1\over 2} C\ell(P) - \tilde\epsilon.$ Clearly,
$\partial \tilde V_{\tilde\epsilon}$ is an annulus and the boundary
components of this annulus are lifts of the meridian of
$V_{\tilde\epsilon}.$ These facts, coupled with the results in [C],
give us very nice bounds on the geometry of $\partial V_{\tilde\epsilon}$.

\medskip
\psfig{file=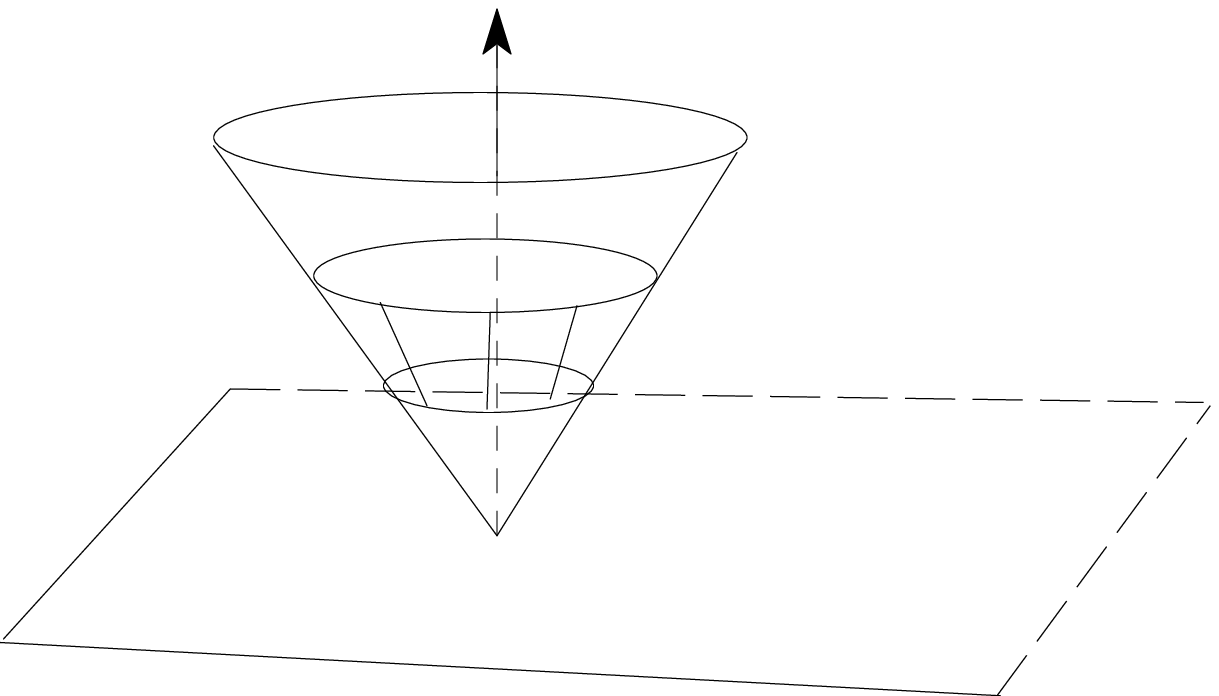,height=2.5in,width=5.0in}
Figure 3.0. {$\tilde V_{\tilde\epsilon}$ in the upper-half space model of ${\Bbb
H}^3$}
\medskip

\proclaim{Lemma 3.3 (Boundary Torus)} Let $M$ be a closed hyperbolic 3-manifold.
Let $P$ be any presentation for $\pi_1(M)$. Suppose $diam(M)\ge C\ell(P)$ where $C>C_1+1$. 
Let $V_{\tilde\epsilon}$ is a deep tube in $M$ given by Proposition 3.2.
The following statements then hold for $\partial V_{\tilde\epsilon}$:
\newline \noindent (i) The induced metric on $\partial
V_{\tilde\epsilon}$ is Euclidean.
\newline \noindent (ii) For every $x\in\partial V_{\tilde\epsilon}$, we 
have $inj(x)\ge \tilde\epsilon$.
\newline \noindent (iii) A loop $\alpha$ on
$\partial V_{\tilde\epsilon}$ homologous to a power of the meridian has\newline
\noindent $length(\alpha) \ge 2\pi sinh({1\over 2} C\ell(P))$.
\endproclaim
\noindent {\bf Proof} Certainly $\partial V_{\tilde\epsilon}$ is
topologically a torus. Thus, (i) will be established if we show that for
each point $x$ in $\partial V_{\tilde\epsilon}$, the curvature at $x$ is
0. We note that $\partial V_{\tilde\epsilon}$ is homogeneous. To see
this, let $x$ and $y$ be distinct points in $\partial
V_{\tilde\epsilon}$. Choose lifts $\tilde x$ and $\tilde y$ of these
points to $\partial \tilde V_{\tilde\epsilon}$. There is a rotation and
translation about the z-axis in ${\Bbb H}^3$ taking $\tilde x$ to
$\tilde y$. Since this commutes with the covering isometry
$\tilde\gamma$, it gives a well-defined isometry of $\partial
V_{\tilde\epsilon}$ which takes $x$ to $y$. So it suffices to prove that
the curvature is 0 at one point in $\partial V_{\tilde\epsilon}$.
Clearly, the curvature of $\partial V_{\tilde\epsilon}$ is not positive
since otherwise the Gauss-Bonnet Theorem implies that $\partial
V_{\tilde\epsilon}$ is a sphere. On the other hand, since $\rho:{\Bbb
H}^3\longrightarrow M$ is a local isometry, the curvature at a point
$\tilde x$ in $\partial \tilde V_{\tilde\epsilon}$ must equal the
curvature at $x=\rho(\tilde x)$ in $\partial V_{\tilde\epsilon}$. But
$\tilde\partial V_{\tilde\epsilon}$ is cylinder with many distinct
simple closed geodesics; in fact, every horosphere centered at $\infty$
intersects $\tilde\partial V_{\tilde\epsilon}$ in a simple closed
geodesic. A second application of Gauss-Bonnet shows this fact implies
the curvature of $\tilde \partial V_{\tilde\epsilon}$ cannot be
negative. This completes the proof since the cylinder itself is
homogeneous. 
Statement (ii) follows at once from our previous discussion. For (iii),
suppose $\alpha$ is a power of the meridian. Since the metric on 
$\partial V_{\tilde\epsilon}$ is Euclidean, this loop is at least as long 
as a geodesic representative of the meridian. A geodesic meridian lifts
to a ``horizontal'' circle; that is, the intersection of a
horosphere centered at $\infty$ with $\tilde\partial
V_{\tilde\epsilon}$. This has length 
$2\pi sinh({1\over 2} C\ell(P))$. $\square$

\medskip

\proclaim {Proposition 3.4} Let $B({\tilde\epsilon\over 2})$ denote the
standard ball of radius ${\tilde\epsilon\over 2}$ in ${\Bbb H}^3$. There 
is a constant $C_2>C_1 + 1$ such that if $M$ is a closed hyperbolic 3-manifold, $diam(M)\ge C_2\ell(P)$, 
and $V_{\tilde\epsilon}$ is a deep tube in $M$, then 
${1\over 2} Vol(B({\tilde\epsilon\over 2})) \le 
Area(\partial V_{\tilde\epsilon})\le Vol(M)$.  Thus, 
$Area(\partial V_{\tilde\epsilon}) \le 2\pi\ell(P).$ 
\endproclaim

\noindent {\bf Proof} A simple calculation shows that for a hyperbolic solid torus of 
radius $L$,

$${Area(\partial V_{\tilde\epsilon})\over Vol(V_{\tilde\epsilon})} = 2{cosh(L)\over sinh(L)}$$

\noindent Thus, it suffices to find bounds on 
$Vol(V_{\tilde\epsilon})$. Choosing a point $x$ on $\partial\tilde 
V_{\tilde\epsilon}$, we have that $inj(x) = \tilde\epsilon$.  Hence, the 
restriction of the universal covering map $\rho:{\Bbb H}^3 
\longrightarrow M$ to $B(x,{\tilde\epsilon \over 2})$ is an embedding.  
If $C_2>C_1 + 1$ is chosen to be sufficiently large, then the tube radius is very 
large. Hence, in our setting $B(x,{\tilde\epsilon\over 2})\cap \partial V_{\tilde\epsilon}$ is
approximately a hemisphere centered at $x$.  This means that for a small
$\delta^\prime$ which depends only on $C_2$,

$$Vol(V_{\tilde\epsilon})\ge {1\over 2}Vol(B(x,{\tilde\epsilon\over
2}))-\delta^\prime$$

\noindent where $\delta^\prime\rightarrow 0$ as $C_2\rightarrow\infty$.
The above shows that we may assume 

$${1\over 2}Vol(B(x,{\tilde\epsilon\over 2}))-\delta^\prime\ge{1\over 4}
Vol(B(x,{\tilde\epsilon\over 2}))$$

\noindent whence 

$$Area(\partial V_{\tilde\epsilon}) \ge 2({1\over 2}Vol(B(x,{\tilde\epsilon
\over 2})) -\delta^\prime) + \delta \ge 2({1\over 4}Vol(B(x,{\tilde\epsilon
\over 2})).$$

\noindent The right hand inequality of then follows since 
$$Area(\partial V_{\tilde\epsilon})\le
2vol(V_{\tilde\epsilon})<2\pi\ell(P).$$ $\square$
\medskip

Now suppose $X$ and $Y$ are closed loops on a Euclidean torus.  
We write $\Delta(X,Y)$ for the algebraic intersection number of 
$X$ and $Y$.
Note that for a pair closed geodesics on a torus, the algebraic and 
geometric intersection numbers are equal.  We prove a simple result
about interesctions of closed loops on a Euclidean torus with 
"bounded geometry" that is quite useful in our setting.

\proclaim {Proposition 3.5} Suppose $T$ is a Euclidean Torus with area
$A$ and injectivity radius $R$.  There is a basis 
$\lbrace [X],[Y]\rbrace$
for $\pi_1(T)$ with $X$ and $Y$ closed geodesics such that:
\newline\noindent (i) $\Delta([X],[Y]) = 1$,
\newline\noindent (ii) $\ell(X) = R$,
\newline\noindent (iii) $\ell(Y) \le {2A\over\sqrt{3}R}$.
\endproclaim

\noindent {\bf Proof} Fix a universal cover $\rho:{\Bbb
R^2}\longrightarrow T$ and choose a point $x$ on $T$.
Since $inj(T) = R$, there is a closed geodesic
loop $X$ based at $x$ with length $R$.  The lifts of $X$
thus form a set of parallel lines in the plane, which for simplicity we 
may assume are horizontal.  For each lift $\tilde X$, the set 
$\tilde X \cap \rho^{-1}(x)$ is a
collection of points such that each adjacent pair is separated by a
distance of $R$. A straightforward 
trigonometry calculation shows that each adjactent pair of lifts, say  
$\tilde X$ and 
$\tilde L$, must have the property that $d(\tilde X,\tilde L)$ 
is at least ${\sqrt{3}\over 2} \epsilon$ as otherwise we can join a pair 
of points in $\rho^{-1}(x)$ by a path of length less than $R$.
This would produce an essential loop in $T$ of length less than
$R$ which is impossible.  On the other hand, it is easy to 
construct a
fundamental domain for the action of $\pi_1(T)$ on ${\Bbb R^2}$
as follows:  choose a point $\tilde x_1$ on $\tilde X\cap
\rho^{-1}(x)$.  There is a point $\tilde y$ on $\tilde L \cap
\rho^{-1}(x)$ with the property that 
$$d(\tilde x_1,\tilde y) = \text{min} \lbrace d(\tilde x_1,z)\text{ }|
\text{ } z\in \rho^{-1}(x) - \tilde x_1)\rbrace.$$
By joining $\tilde x_1$ and $\tilde y$ with a straight segment and 
then constructing a parallel segment between adjacent points on 
$\tilde X$ and $\tilde L$ respectively, we obtain a parallelogram which 
is a fundamental domain for $\pi_1(T)$.  Also, the segment joining $x_1$ 
and $y$ projects to an essential loop $Y$ on $T$ which by construction 
has $\Delta([X],[Y])=1$.  By a short calculation, we conclude 
that $\ell(Y) \le {2A\over\sqrt{3}R}$, since the area of
this parallelogram is exactly $A$. $\square$

\medskip 
\proclaim{Definition 3.6} Suppose $T$ is a Euclidean Torus with area
$A$ and injectivity radius $R$. A basis $\lbrace [X],[Y]\rbrace$ for
$\pi_1(T)$ that satisfies the conditions of Proposition 3.5 is
a {\bf short basis} for $\pi_1(T)$.
\endproclaim
\medskip

\proclaim{Proposition 3.7} Suppose $T$ is a Euclidean Torus with injectivity radius $R$. 
Let $L$ be an essential loop on $T$. If $\lbrace [X],[Y] \rbrace$ is a short basis for
$\pi_1(T)$, then $[L] = a[X] + b[Y]$ where $\text{max } \lbrace
|a|,|b|\rbrace \le {2\text{ } length(L)\over \sqrt{3} R}$.
\endproclaim

\noindent {\bf Proof} We can tile the plane with copies of the 
parallelogram fundamental domain for $T$ as constructed in the 
proof of Proposition 3.5.  Since the distance between parallel sides of
a given parallelogram is at least ${\sqrt{3}\over 2} S$, the lift of a 
geodesic representative of $[L]$ hits at most 
$2\text{ }length(L)\over \sqrt{3} R$ parallelograms. In particular, 
$\Delta([L],[X])$ and $\Delta([L],[Y])$ are both bounded above by
$2\text{ } length(L)\over\sqrt{3} R$. $\square$ 

We can apply the above results to our setting at once since the 
injectivity radius on $\partial V_{\tilde\epsilon}$ is 
approximately the Margulis constant $\tilde \epsilon$.  In fact the
injectivity radius on the boundary torus is slightly larger.  To see this, 
lift a short essential loop on $\partial V_{\tilde\epsilon}$ to ${\Bbb H_3}$. 
This gives a path between distinct points on 
$\tilde \partial V_{\tilde\epsilon}$.  The geodesic path in ${\Bbb H^3}$
joining these points has length at least $2\tilde\epsilon$, so that the
lifted path is slightly larger.

\proclaim {Corollary 3.8} Let $V$ denote a margulis tube of $M$.  
There exists a short basis $\lbrace [X],[Y] \rbrace$ for $\pi_1(\partial V_{\tilde\epsilon})$ 
such that if $L$ is any essential loop on $\partial V$, then $[L] = a[X] + b[Y]$ where 
$\text{max }\lbrace |a|,|b|\rbrace \le {2 length(L)\over \sqrt{3} \tilde\epsilon}$.
\endproclaim

Suppose that $f:K\longrightarrow M$ is a triangular complex. Lemma 2.2 
(Invariant Intersection) shows that $f(K)$ must meet $V$.  This gives us
two important in subsets of $K$. Since $f(K)$ must meet the core of $V_{\tilde\epsilon}$, the set 
$f^{-1}(V_{\tilde\epsilon})$ is mapped very far into the deep tube.  Define the {\bf pull-back
boundary} of $f^{-1}(V_{\tilde\epsilon})$ to be the set 
$f^{-1}(\partial V_{\tilde\epsilon})$. It follows that
the pull-back boundary is a subset of the topological boundary of $f^{-1}(V_{\tilde\epsilon})$.
However, we note that in general the pull-back boundary does not equal
the topological boundary.
Lemma 3.3 (Boundary Torus) allows us to completely describe 
these sets geometrically.  We obtain three types of intersection set as shown in Figure 3.1.  
We call these intersection sets {\bf 0-handles}, {\bf 1-handles}, and 
{\bf monkey-handles} respectively.  Notice that a $monkey$-handle
appears as a disc with six edges, three of which join adjactent edges of
the triangle.  A given triangle can contain at most one $monkey$-handle.

\medskip
\psfig{file=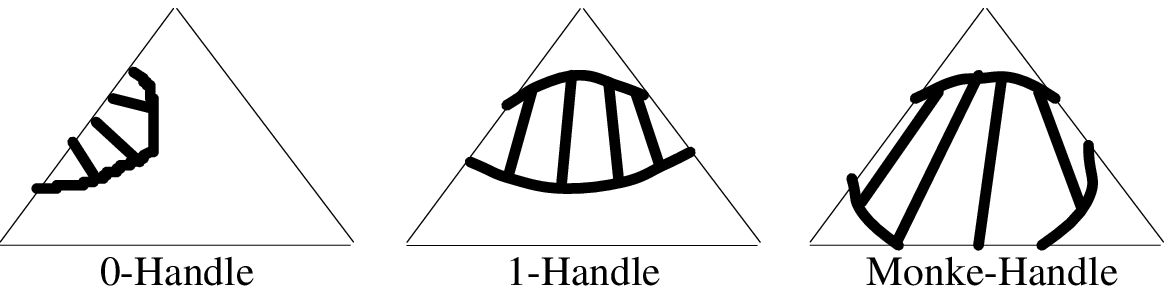,height=1.25in,width=5.0in}
Figure 3.1.
\medskip

\proclaim{Lemma 3.9} Let $f:K\longrightarrow M$ be a triangular
complex. Let $V$ be any Margulis tube in $M$. 
Then:\newline \noindent (i) $f^{-1}(V_{\tilde\epsilon})$ is a finite union of $0$-handles,
$1$-handles, and $monkey$-handles, together with a collection of discs contained in the 
interiors of triangles of $K$.\newline
\noindent (ii) The pull-back boundary is a 
(possibly not connected) graph.\newline \noindent (iii) The number of 
$monkey$-handles in $f^{-1}(V_{\tilde\epsilon})$ is at most the number of triangles in $K$. 
\newline
\noindent (iv) Removal of the interiors of the $monkey$-handles from 
$f^{-1}(V_{\tilde\epsilon})$ yields an I-bundle over a graph. \endproclaim 

\noindent {\bf Proof} We may assume that the single vertex of $K$ is
mapped by $f$ into the thick part of $M$.  Since every relation in $\pi_1(M)$ corresponds to 
a disc ${\Bbb D}$ in $K$ which lifts to a hyperbolic triangle in 
${\Bbb H}^3$, we may view ${\Bbb D}$ as a hyperbolic triangle with 
vertices identified to a single point. In particular, we can understand 
$K_v$ by looking at the intersection of each lifted triangle in 
${\Bbb H^3}$ with
$\tilde{\partial V_{\tilde\epsilon}}$. By an arbitraritly small adjustment,
we can assume that the intersection of a hyperbolic triangle with the 
conical boundary,$\partial V_{\tilde\epsilon}$, is transverse. Both $\tilde V_{\tilde\epsilon}$
and the triangle are convex.  This implies their intersection is convex, so the intersection
must be a disc.  Also, each triangle edge is convex, so the intersection of each edge with
this disc is either empty or an interval.  We readily obtain the three types
pictured in Figure 3.1. To see that the union is finite, consider a
fundamental domain for $V_{\tilde\epsilon}$. A fundamental domain for $K$ 
in ${\Bbb H}^3$ is a finite union of (compact) hyperbolic triangles.  The
action of $\pi_1(M)$ on ${\Bbb H}^3$ is properly discontinuous.
Therefore, only finitely many translates under $\pi_1(M)$ of a given
triangle in $K$ hit the fundamental domian for $V_{\tilde\epsilon}$
This proves (i) and (ii). Part (iii) is obvious. Also, (iv) follows 
since the complement of the interiors of the 
$monkey$-handles in $K_v$ is built by glueing 0-handles and 1-handles
(which are I-bundles over an interval) along interval fibres. $\square$.

\medskip 
\noindent {\bf Remark} In our subsequent arguments, we do not need to consider discs in $f^{-1}(V_{\tilde\epsilon})$ 
which are contained in the interior of a triangle of $K$.  Thus, we let 
$K_V=f^{-1}(V_{\tilde\epsilon}) - \lbrace\text{interior discs}\rbrace$ and define $\partial_f K_V$ to be the
corresponding pull-back boundary. 

\medskip
\psfig{file=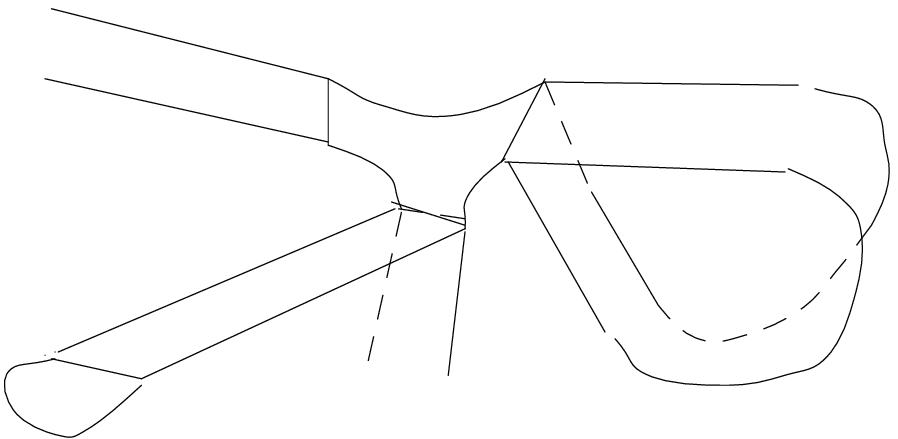,height=2.25in,width=5.0in}
Figure 3.2. A piece of $K_V$. Note that removing the monkey-handle interiors gives an
I-Bundle over a graph.
\medskip

\noindent We use Proposition 2.1 (Bounded Area) to show that $K_V$ has bounded
geometry.  To explain this, let $\lbrace \partial V_s \rbrace_{s=0}^1$ be 
a parameterization by distance $s$ of the parallel concentric tori in
$\partial V_{\tilde\epsilon}$; that is, $\partial V_s$ is the parallel
torus of distance $s$ from $\partial V_{\tilde\epsilon}$. Let $\Gamma
(s) = K_v \cap f^{-1}(\partial V_s)$. We can show that there is an $s$
such that $0 < s < 1$ and $length(\Gamma (s))$ is bounded by presentation
length. Thus, by replacing $V_{\tilde\epsilon}$ with $V_s$ if necessary, we may assume that
the length of $\partial_f K_v$ is also bounded in terms of presentation length. In 
section 4, we shall be able to prove results bounding the homology of 
$K_v$ using this fact. 

\medskip
\psfig{file=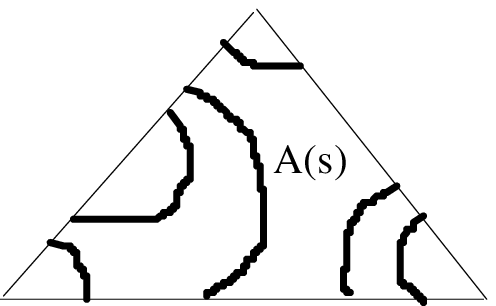,height=1.5in,width=2.0in}
Figure 3.3.
\medskip

\proclaim{Lemma 3.10 (Short Boundary)} With the notation as above, there
is a universal constant $C_3 > C_2$ such that $diam(M) > C_3\ell(P)$
implies there exists $s_0$ with $0 < s_0 < 1$ and $length(\Gamma(s_0)) <
area(K)$. Thus, $length(\Gamma(s_0)) < 2\pi\ell(P)$. \endproclaim 

\noindent {\bf Proof} It is sufficient to show this for a single triangle in $K$.  So let
$T$ be a triangle in $K$ and consider $A(s) = \Gamma(s)\cap T$.  Then $A(s)$ is a collection of
embedded arcs in $T$. The endpoints of each arc are contained in edges of $T$.
At most three of these are tangent to an edge of $T$; all others are properly embedded.  
It follows that:

$$\int_0^1 \text{length} (A(s)) ds \le Area (T).$$

\noindent Since $Area(K)$ is the sum of the areas of the triangles in $K$, 
the lemma follows at once. $\square$
\medskip

Given this result, it is best to introduce new notation: let $V$ denote the 
Margulis Tube in $V_{\tilde\epsilon}$ with boundary
$\partial V = \partial V_{s_0}$, where $s_0$ is provided by Lemma
3.5 (Short Boundary).  It is clear from the above proof that we can choose $s_0$ so that $V$ is
properly contained in $V_{\tilde\epsilon}$. This assumption simplifies later exposition since
$V_{\tilde\epsilon}-V$ is homeomorphic to the product of a torus and an interval.
We shall now set $K_V = f^{-1}(V)- \lbrace\text{interior discs}\rbrace$ so that 
$\partial_f K_V < 2\pi\ell(P)$.

\proclaim{Corollary 3.11} Let $M$ be a closed hyperbolic 3-manifold, $P$ a triangular presentation
for $\pi_1(M)$ and $f:K\longrightarrow M$ a corresponding triangular
complex. Assume that $L\subset \partial_f K_V$ is an embedded loop
with $f_*([L])\ne 1$ in $\pi_1(\partial V)$. Then
$diam(M)>C_3\ell(P)$ implies $f_*([L])$ is not a power of the meridian
of $V$. \endproclaim

\noindent {\bf Proof} By conclusion (iii) of Lemma 3.3 (Boundary Torus), 
if $f_*([L])$ is a power of the meridian, then we have 
$length(f(L))\ge 2\pi sinh({1\over 2}C_3\ell(P))$. By our construction, $C_3 > 1$, 
so certainly $2\pi\ell(P)<2\pi sinh({1\over 2}C_3\ell(P))$.  But this implies 
$length (L) = length (f(L))\ge 2\pi\ell(P)$ which contradicts 
Lemma 3.10 (Short Boundary). $\square$ 

\proclaim{Corollary 3.12} A simple loop in $\partial_f K_V$ is
inessential in $K$ if and only if its image under $f$ is inessential in
$\partial V$. \endproclaim

\noindent {\bf Remark} By the above results, $K$ does not contain a disc which is mapped to 
a meridian of a deep tube (since this gives enormous area).  Also, the image of $K$
intersects every essential loop in $V$.  However, as the following
example shows, these facts do not yield a proof.  Using the short basis
for $\pi_1(\partial V)$ provided by Corollary 3.8, sweep out a 2-complex
which hits the core of $V$ essentially by concentrically shrinking the loops $X$ and
$Y$ to the core of $V$.  By construction, this complex has boundary length
$length(X) + length(Y)$ which by Propostions 3.31 and 3.32 is at most 
$\tilde \epsilon + {4\pi\ell(P)\over\sqrt{3}\tilde\epsilon}$.  This means a meridian is too
long to be contained in the boundary of the complex.  On the other hand, this complex has 
bounded area and intersects every essential loop in $V$.  The next two sections deal with 
this complication.

\head 4. Bounds On First Homology 
\endhead
\medskip
In this section, we use the results of Section 3 to establish a bound on
the ``complexity'' of $K_V$. This amounts to
showing that $Rank(H_1(K_V))$ is small compared to diameter. 

\proclaim{Lemma 4.0} Let $G_1,\dots G_r$ denote the components of
$\partial_f K_V$. There exists a component $G_i$ such that the restriction
$f_{G_i}:G_i\longrightarrow \partial V$ has nontrivial induced
homomorphism ${f|_{G_1}}_*:\pi_1(G_i)\longrightarrow \pi_1(\partial
V)$. \endproclaim
\noindent{\bf Proof} If every component of $\partial_f K_V$
has trivial induced homomorphism into $\pi_1(\partial V)$, we may
homotop $f$ on each component $G_i$ so that the
image of $G_i$ under $f$ is a finite collection of points in $\partial
V$. This contradicts Proposition 2.2 since there is a loop in
$\partial V$ that is essential in $M$ and misses each of these points.
$\square$ 

\proclaim{Lemma 4.1} Suppose $f:K\longrightarrow M$, $K$ is a
2-complex, and $f$ is a $\pi_1$- isomorphism. Let ${\Bbb D}$ be an
embedded disc in $K$ with the property that $\partial {\Bbb D}\simeq *$
in $K-int({\Bbb D})$ and $closure(K-{\Bbb D})\cap {\Bbb D} = \partial {\Bbb D}$.
Then the restriction $f|_{K-int({\Bbb D})}:K-int({\Bbb D})\longrightarrow M$ 
is a $\pi_1$-isomprophism.
\endproclaim

\noindent {\bf Proof} Let $i:K-int({\Bbb D})\longrightarrow
K$ and $j:\partial {\Bbb D}\longrightarrow K$ denote the inclusion
mappings. Let $x$ be a basepoint in $K-int({\Bbb D})$ and let
$N(j_*(\pi_1(\partial {\Bbb D},x)))$ denote the smallest normal subgroup of
$\pi_1(K,x)$ containing $j_*(\pi_1(\partial {\Bbb D},x))$. Van Kampen's
Theorem gives that $i_*:\pi_1(K-int({\Bbb D}),x)\longrightarrow
\pi_1(K,x)$ is an epimorphism with kernel $N(j_*(\pi_1(\partial {\Bbb
D},x)))$. The hypothesis gives that $N(j_*(\pi_1(\partial {\Bbb D},x)))
= 1$, so that $i_*$ is an isompophism. This completes the proof since
$f_*$ is an isomorphism and ${f|_{K-int({\Bbb D})}}_* = f_*i_*$
$\square$

\proclaim{Lemma 4.2} Let $K$ be a 2-complex, $M$ a 3-manifold, and
$f:K\longrightarrow M$ a $\pi_1$ isomorphism.  Suppose $\gamma$ is an
inessential loop in $K$ and form the complex $K^\prime = K\cup_\gamma
{\Bbb D}$.  There is a map $f^\prime:K^\prime\longrightarrow M$ that is
a $\pi_1$ isomorphism.
\endproclaim

\noindent {\bf Proof} Since $\gamma$ is inessential, Van Kampen's
theorem gives that $\pi_1(K^\prime)\cong\pi_1(K)$. Also, $f(\gamma)$ is
inessential in $M$, so there is a map $h:{\Bbb D}\longrightarrow M$ with
$h(\partial {\Bbb D})=f(\gamma)$.  Extend $f$ over ${\Bbb D}$ using $h$ and the
proof is complete. $\square$

\proclaim{Lemma 4.3} Suppose $f:{\Bbb D}\longrightarrow
\partial V$. Then there is a map $g:{\Bbb D}\longrightarrow
V_{\tilde\epsilon}$ such that the following conditions hold:\newline
\noindent (i) $g=f$ on $\partial {\Bbb D}$\newline
\noindent (ii) $g\sim f$ rel $\partial {\Bbb D}$\newline
\noindent (iii) $g(int({\Bbb D}))\subset V_{\tilde\epsilon} - V$.
\endproclaim
\noindent {\bf Proof} Push the interior of $f({\Bbb D})$ out into the collar of
$\partial V$ contained in $V_{\tilde\epsilon} - V$. $\square$
\medskip

\noindent We can actually make many such maps 
with that property that any two maps send the interior of the disc to
disjoint sets in $V_{\tilde\epsilon} - V$.  We obtain the following
useful corollary as a consequence:

\proclaim{Corollary 4.4} Suppose $\gamma$ is an inessential loop in
$\partial_f K_V$.  Let $K^\prime = K\cup_\gamma {\Bbb D}$.  Then there is a
$\pi_1$ isomorphic map $f^\prime:K^\prime\longrightarrow M$ such that
$K_V^\prime = K_V$.
\endproclaim

\noindent {\bf Proof} Let $g:{\Bbb D}\longrightarrow V_{\tilde\epsilon}$
denote the map provided by Lemma 4.3.  Then the map $f^\prime = f\cup g$
satisfies our requirements. $\square$
\medskip

Let $K_I$ denote $closure(K_V -
\lbrace\text{interior(monkey-handles)}\rbrace)$. We also define ${\Cal
N}(K_I)$ by taking $K_I$ together with a small product neighborhood of
its boundary in $K$. We note in passing that this definition ensures
that the pull-back boundary of this neighborhood $\partial_f {\Cal N}(K_I)$ maps 
into $V_{\tilde\epsilon} - V$. By Lemma 3.9, $K_I$ is an I-bundle over a graph.  
We wish to show
that we can modify $K$ so that $rank(H_1(K_I))$ is bounded by $\ell(P)$.

To do this, we are going to define a surgery procedure made possible by
the above two lemmas.  We are motivated by two of our constraints: $length (\partial_f K_V)$ is 
bounded above and $inj(\partial_f K_V)$ is bounded below.  Thus, if the rank of $H_1(K_I)$ is 
large, then there are many short inessential loops in $\partial_f K_V$.
The point of our procedure will be to ``snap'' all of the
annuli and mobius bands in $K_I$ that are inessential in $K$ by removing a single
1-handle from each. To preserve the $\pi_1$ isomorphism, we must make
sure that the resulting new loop remains inessential in the new
complex $K-\lbrace 1-handle\rbrace$. We treat the annulus case first.
We can attach a pair of discs to the boundary of the original
inessential annulus in ${\Cal N}(K_I)$. By Lemma 4.3, we can ensure
that the interiors of these discs miss $V$, so that $K_V$
remains unchanged; that is, we glue a disc onto each boundary component
to produce a 2-sphere. Removal of the 1-handle thus gives an
inessential loop. The effect is that $K_I$ appears to lose a
1-handle. 

\medskip
\psfig{file=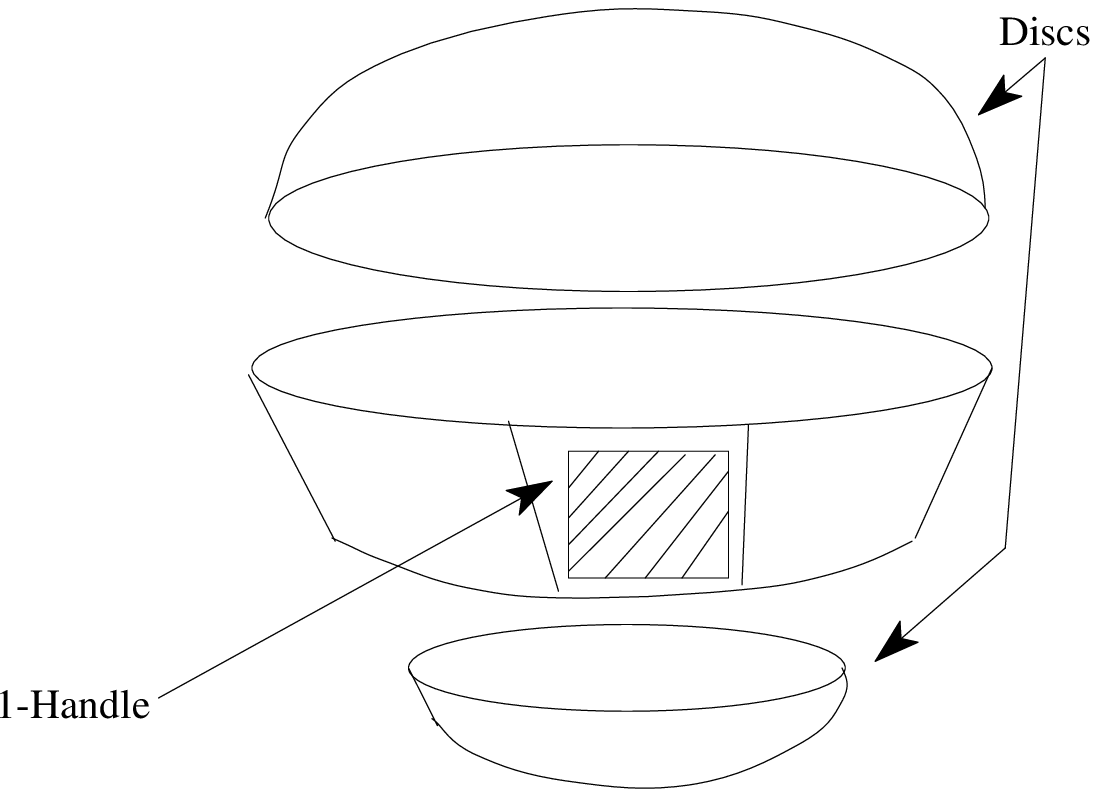,height=3.0in,width=5.0in}
Figure 4.0.
\medskip

\proclaim{Lemma 4.5 (Inessential Annulus)} Suppose $A$ is an
embedded annulus in $K_I$ with $\partial A\subset\partial_f K_I$. 
Assume that $\partial A\simeq *$ in $K$. Let
$H$ denote any one of the 1-$handles$ contained in $A$. Define $H^+$ to
be the union of $H$ and a small neighborhood of $\partial_f K_I \cap H$ in
the triangle of $K$ containing $H$. There is a 2-complex $K^\prime$ with
$K-int(H^+)\subset K^\prime$ and a map $f^\prime:K^\prime-int(H^+)\longrightarrow
M$ is a $\pi_1$ isomorphism. Moreover, $K^\prime_V = 
K_V - H$.
\endproclaim 

\noindent {\bf Proof} By Corollary 3.12, both components of $\partial A$
map to loops that are inessential in $\partial V$. Using
Corollary 4.4 twice, we may attach a disc to each boundary componet of
$A\cup H^+$ to build the $\pi_1$ isomorphic 2-complex
$f^\prime:K^\prime\longrightarrow M$. The union of these two discs with
$A\cup H^+$ is an embedded $S^2$ in $K^\prime$. Therefore, $\partial
H^+$ is inessential in $K^\prime - int(H^+)$. By Lemma 4.1,
$f^\prime:K^\prime-int(H^+)\longrightarrow M$ is a $\pi_1$ isomorphism.
It is then obvious that $K^\prime_I = K_I - H$. $\square$ 
\medskip

\noindent {\bf Remark} Since we define the pull-back boundary as 
$\partial_f K_V = f^{-1}(\partial V)$, application of 
Lemma 4.5 (Inessential Annulus) does not increase the length of $\partial_f K_V$.

\medskip

The mobius band case is more complicated, but can be approached in a
similar way.  We say that an embedded mobius band $B$ in 
$K_V$ with $\partial B\subset\partial_f K_V$ is {\bf inessential} if $\partial B$ is inessential 
in $K$. This terminology is justified by the following propostion:

\proclaim{Proposition 4.6} Suppose $B$ is an inessential mobius band in $K_V$. Let $\gamma$ 
denote the core curve of $B$. Then $\gamma\simeq *$ in $K$. 
\endproclaim

\noindent {\bf Proof} Since $\gamma^2$ is homotopic to $\partial B$, 
$[\gamma]^2 = 0$ in $\pi_1(M)$.  Since $M$ is a closed hyperbolic 3-manifold, $\pi_1 M$ is 
torsion free. The map $f:K\longrightarrow M$ is a $\pi_1$-isomorphism, so $\gamma$ is 
contractible in $K$. $\square$

\medskip

\noindent Suppose $B$ is an inessential mobius band in $K_I$.  We glue two discs onto $B$ to 
obtain $f:K^\prime\longrightarrow M$ such that the first disc contributes a 0-handle to 
$K^\prime_I$ and the second disc is mapped into $V_{\tilde\epsilon}-V$ (and thus is disjoint
from $K^\prime_I$).  This is done in such a way that removing a 1-handle from $B$ in $K^\prime$
does not change $\pi_1 K^\prime$. Since we plan to repeat this procedure for each inessential 
mobius band, we must ensure that the new attached 0-handle does not increase
$length (\partial_f K_I)$ too much. Let $B$ denote an inessential mobius band 
in ${\Cal N}(K_I)$. Thus, every closed loop in $B$ is null 
homotopic in $K$. 

\proclaim {Definition 4.7} A {\bf good core} $\gamma$ for $B$ is an
embedded loop such that the following conditions hold: 
\newline\noindent (i) $[\gamma]$ generates $\pi_1(B)$, 
\newline\noindent (ii) $\gamma$ is the union of a pair of arcs 
denoted ${\Cal A} \cup {\Cal L}$ where ${\Cal A}$ is an arc on 
$\partial {\Cal N}(K_I)$, ${\Cal L}$ is a properly embedded geodesic arc 
in $B$, and ${\Cal A} \cap {\Cal L}$ is a pair of distinct points
(see Figure 4.1), 
\newline\noindent (iii) $f(\gamma)$ is embedded in $M$. 
\endproclaim

\noindent We use condition (iii) since we will attach a disc to the good core: we must control
the resulting intersection of the image of this disc with our deep tube.

\medskip
\psfig{file=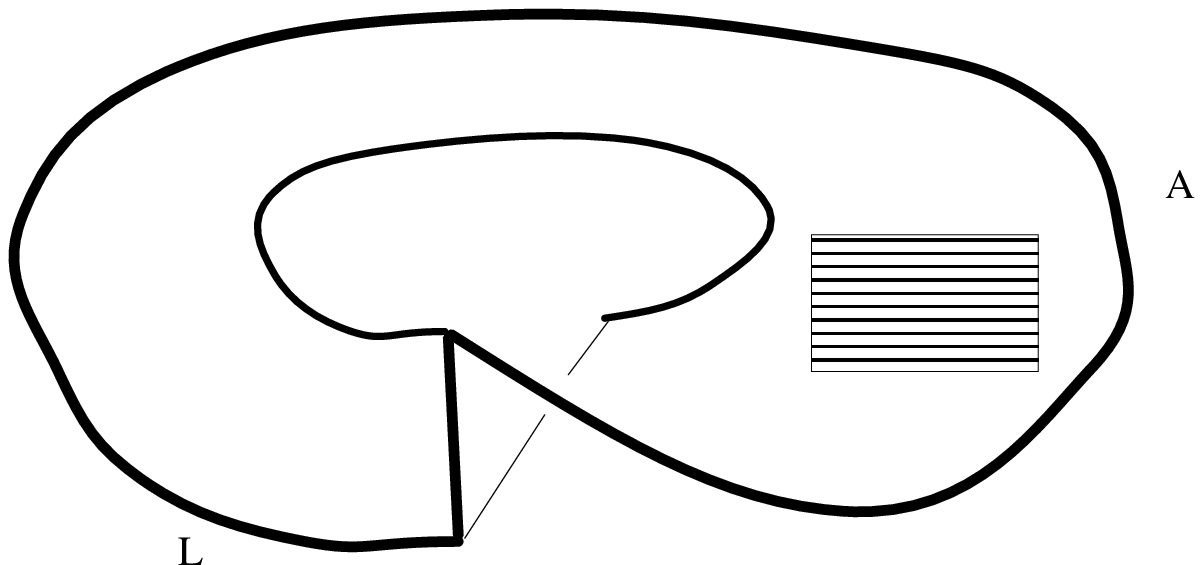,height=2.0in,width=5.0in}
Figure 4.1.
\medskip

Since $f(\gamma)\simeq *$, there exists a lift $\tilde\gamma$ to ${\Bbb H}^3$
which is an embedding into $\tilde V_{\tilde\epsilon}$.  Notice that
$\tilde V_{\tilde\epsilon}$ and $\tilde V$ are a pair of 
concentric cones about the z-axis in ${\Bbb H}^3$. ${\Cal L}$ lifts to 
an arc $\tilde {\Cal L}$ which is properly embedded in $V$
while ${\Cal A}$ lifts to an arc which is embedded in a concentric cone 
about the z-axis in $\tilde V_{\tilde\epsilon} - int(V)$.  
The image $\tilde {f(\gamma)}$ thus bounds an embedded disc in ${\Bbb H}^3$ 
which intersects $V$ in a single 0-handle.  We construct this 
disc as follows: $\tilde {\Cal A}$ intersects $\partial \tilde V$
in a pair of points.  Choose the shortest path ${\Cal P}$ on $\partial 
\tilde V$ joining these two points.  Notice that by construction
$length({\Cal P}) < length(\tilde {\Cal A}) < length (\partial B)$.  

Now $\tilde {\Cal A}\cup {\Cal P}$ is a circle which bounds an embedded 
disc ${\Bbb D}_1$.  Likewise, ${\Cal P}\cup {\Cal L}$ bounds an embedded
disc ${\Bbb D}_2$.  Put ${\Bbb D} = {\Bbb D}_1 \cup {\Bbb D}_2$.
Then, by Lemma 4.2, we may attach ${\Bbb D}$ to $\gamma$ to obtain a 
2-complex $f^\prime:K^\prime \longrightarrow M$ with $f^\prime_*$ a 
$\pi_1$ isomorphism and $K^\prime_I = K_I \cup \lbrace$0 -$handle\rbrace$.
Furthermore, $length (\partial_f K^\prime_I) < length (\partial_f K_I) + length(B)$. 
Our surgery argument is then completed by the next proposition.
\medskip 
\proclaim{Proposition 4.8} Suppose $B$ is a mobius band.  Let ${\Bbb
D}_X$ be an embedded disc in $B$ with boundary $X$.
Assume that $interior{\Bbb D}_X$ is disjoint from $\partial B$ and from
a good core $\gamma$.  Form the space ${\Cal X}$ by attaching one disc 
${\Bbb D}_\gamma$ to the core curve and another disc ${\Bbb D}_{\partial
B}$ to the boundary of the mobius band.  Then ${\Cal X} - interior({\Bbb D}_X)$
is simply connected.
\endproclaim
\medskip
\noindent {\bf Proof} There is a retraction of $B - interior({\Bbb D}_X)$
onto a closed trivalent graph $G$ which has fundamental group ${\Bbb
Z}*{\Bbb Z}$.  This retraction maps $\partial B$ and $\gamma$ to a pair
of generators for $\pi_1(G)$.  Therefore, the space ${\Cal X} - 
interior({\Bbb D}_X)$ is homeomorphic to $G$ with a pair of discs glued 
onto a generating set for
$\pi_1(G)$ so that $\pi_1({\Cal X} - interior({\Bbb D}_X))$ is trivial.
$\square$
\medskip

\proclaim{Lemma 4.9 (Inessential Mobius Band)} Suppose $B$ is an embedded inessential mobius
band  in $K_I$. Let $H$ denote any one of the 1-$handles$ contained in $B$. Let $H^+$ denote the
the union of $H$ and a small neighborhood of $\partial_f K_I \cap H$ in
the triangle of $K$ containing $H$. There is a 2-complex $K^\prime$ with
$K-int(H^+)\subset K^\prime$ and a map $f^\prime:K^\prime-int(H^+)\longrightarrow
M$ which is a $\pi_1$ isomorphism. The map $f^\prime$ agrees with $f$ on $K-int(H^+)$.
Moreover, the I-bundle part of $K^\prime_V$ is $K^\prime_I = K_I - H$ together 
with a single 0-handle attached along a proper geodesic edge of $K_I$.
Furthermore, $length (\partial_f K^\prime_V) \le length (\partial 
K_V) + length (\partial B)$.
\endproclaim 

\noindent {\bf Proof} Lemmas 4.1 and 4.2 imply that if we attach
discs to $\partial B$ and a good core $\gamma$ which misses $H$, then 
removal of the interior of $H$ does not change the fundamental group.
The existence of $K^\prime$ as stated follows from the discussion above.
$\square$
\medskip
Noting that, by Lemma 3.9 (i), $K_I$ is composed of finitely many
handles we have:
\medskip
\proclaim {Lemma 4.10} Suppose $M$ is a closed hyperbolic 3-manifold,
$P$ is a triangular presentation of $\pi_1(M)$, and $diam(M) > C_3\ell(P)$.  
Let $V$ be a deep tube in $M$  There is a 2-complex $K^\prime$ and a $\pi_1$-isomorphic map 
$f^\prime:K^\prime\longrightarrow M$ constructed by surgery on 
$K$ such that the subcomplex $\partial_f K^\prime_V = f^{\prime -1}(V)$
has the following property. If $K^\prime_I$ the I-bundle part $K^\prime_V$, 
then every properly embedded annulus and mobius band in 
$K^\prime_I$ is essential in $K^\prime$ and $length(\partial_f K^\prime_V)
\le 2\text{ } length(\partial_f K_V)$.
\endproclaim
\medskip
\noindent {\bf Proof} Apply Lemma 4.5 so that every embedded annulus in $K_I$ is essential. 
Next, apply Lemma 4.9 repeatedly to construct $K^\prime$.  For each mobius 
band $B$, the use of Lemma 4.9 adds at most $length(\partial B)$ to 
the length  of the boundary of the complex.  Notice that since Lemma 4.5
has already cut every inessential annulus in $K_I$, there does not
exist a pair of inessential mobius bands which share a common 1-handle. For, 
if such a pair exists one sees an inessential annulus by removing the common 1-handle.
Therefore, inessential mobius bands are disjoint.  This means that 
if we sum the lengths of all the inessential mobius bands in $K_I$, 
Lemma 3.10 (Short Boundary) implies:
$$length(\partial_f K_V^\prime) < length(\partial_f K_V) + \sum_{B\subset K_I}
length(\partial B) < 2\text{ }length(\partial(K_V))$$
\noindent $\square$
\medskip
The principal use of this fact is to get a bound on the first homology
rank of the I-bundle portion of the 2-complex. We do this by taking
stock of our various results bounding length. First,
the Margulis Lemma implies that essential loops on $\partial_f K_V$
always have length at least the Margulis Constant $\tilde \epsilon$.  
Second, Lemma 4.6 tells us that we may assume every annulus and mobius 
band in $K_I$ is essential.  So each has length at least $\tilde
\epsilon$.  Moreover, we are guaranteed that $length (\partial_f K^\prime)
\le 2 length (\partial_f K_V)\le 4\pi\ell(P)$.  In our situation, the number of annuli and
mobius bands gives a bound on $Rank(H_1(K_I))$.  Thus, we have the
following:

\proclaim {Lemma 4.11} Let $G$ be a finite metric graph with the following properties:\newline
\noindent (i) $length(G) < N$ \newline
\noindent (ii) Every simple closed curve in $G$ has length at least
$\epsilon$. \newline
\noindent Then $Rank(H_1(G)) \le {32N^2\over\epsilon^2}$.
\endproclaim

\noindent {\bf Proof} Let $T$ be a maximal tree of $G$.  We wish to bound the number of edges
in $G-T$ since this bounds the first homology rank.  Now, given any integer $m>0$, condition (i)
implies $G-T$ contains at most $m$ edges of length at least $N\over m$.  So if we choose 
$m={5N\over\epsilon}$, $G-T$ contains at most $5N\over\epsilon$ edges of length at least
$\epsilon\over 5$.  Hence, let ${\Cal L}$ denote the collection of edges in $G-T$ of length 
less than $\epsilon\over 5$.  If ${\Cal L}$ is empty or contains only one edge, the proof is 
complete.  Otherwise, given any pair of edges $A$ and $B$ in $G-T$, one vertex of $A$
has distance at least 
$${1\over 2}(\epsilon - length(A) - length(B)) \ge 
{1\over 2}(\epsilon - 2{\epsilon\over 5}) = {3\epsilon\over 10}$$

\noindent from one vertex of $B$.  This follows since otherwise we can contsruct a simple closed 
curve of length less than $\epsilon$ which contradicts condition (ii). 
Choose a maximal set of vertices $S$ in $T$ such that for every pair of
vertices $v$ and $w$ in $S$, we have $d(v,w) \ge {2\epsilon\over 5}$.
We may center a family of pairwise disjoint balls of radius 
$\epsilon\over 5$ at the vertices in $S$.  Now notice that each
such ball must contain an edge path of length at least
$\epsilon\over 5$.  Since the balls are pairwise disjoint,
there are at most $m$ such balls. Thus, none of the edges in ${\Cal L}$ has both endpoints in a
single ball as this would also violate condition (ii).  On the
other hand, if two of the edges of ${\Cal L}$ connect the same pair of balls, 
there is a simple closed curve of length less than $\epsilon$
which also contradicts condition (ii). Therefore, there are at most 

$$\pmatrix m\cr 2 \endpmatrix = {1\over 2}m(m-1) \le {1\over 2}({8N\over\epsilon})^2$$

\noindent edges in $G-T$ which completes the proof. $\square$
\medskip

\noindent It now follows easily that the first homology of 
$K^\prime_V$ is similarly bounded. 

\proclaim {Theorem 4.12 (Bounded Homology)} Suppose $M$ is a closed hyperbolic 3-manifold,
$P$ is a triangular presentation of $\pi_1(M)$, and $diam(M) > C_3\ell(P)$.  
Let $V$ be a deep tube in $M$. There is a 2-complex $K$ and a $\pi_1$-isomorphic map 
$f:K\longrightarrow M$ such that the subcomplex $\partial_f K_V = f^{-1}(V)$  
has $Rank(H_1(K^\prime_V)) \le B_1 \ell(P)^2$ where $B_1 = ({128\pi^2\over\tilde\epsilon^2} + 3)$ and 
$length(\partial_f K_V)\le 4\pi\ell(P)$.
\endproclaim

\noindent{\bf Proof}  Notice that $K_I^\prime$ can be retracted onto a spine $G$ which is
a graph with the following properties: $length(G) < length(\partial_f K_I^\prime) < 4\pi\ell(P)$ and
every simple closed curve in $G$ has length at least $\tilde \epsilon$.  Hence Lemma 4.11 gives
that $Rank(H_1(K_I^\prime)) \le {32(4\pi\ell(P)^2)\over \tilde \epsilon^2}$.
Now $K^\prime_V$ is built by attaching 
at most $\ell(P)$ $monkey$-handles to $K^\prime_I$. 
Attaching each $monkey$-handle corresponds to attaching a trivalent vertex to the spine; 
so we can bound the homolgy of $K^\prime_V$ by watching what 
happens as we attach tiny neighborhoods of trivalent vertices to the 
spine of $K_I$.  It is an easy argument that each trivalent vertex 
addition increases the first homology rank by at most 3. Thus, the 
lemma is proved.  $\square$

While the above theorem tells us that we $K^\prime_V$ has bounded first
Betti number, we have done nothing to bound the complexity of $\partial
K^\prime_V$. We have no bound on the number of 0-handles in
$K^\prime_V$. Indeed, we have added 0-handles during our surgery
procedure. Each 0-handle contributes one edge to $\partial_f K^\prime_V$,
so it is possible that $\partial_f K^\prime_V$ has a large number of short
edges. For our subsequent arguments, it is also necessary to bound the
number of edges in $\partial_f K^\prime_V$. Hence, we introduce a final
surgery procedure to remove all but a bounded number of 0-handles. As
shown in our previous discussion, a zero 0-handle in $K^\prime_V$ is a
disc with contained in a single triangle of $K^\prime$. The boundary of
the 0-handle consists of two edges. One edge is a subarc of an edge of
the triangle. The other edge is contained in $\partial_f K^\prime_V$.
Consider the special case in which two 0-handles are joined together by
two triangles of $K^\prime$ glued along a common edge as shown in Figure 4.2. This
gives a disc ${\Bbb D}$ in $K^\prime_V$ with $\partial {\Bbb D}$ an
embedded loop in $\partial_f K^\prime_V$. As in the proof of Lemma 4.5
(Inessential Annulus), using Corollary 4.4, attach a disc to $\partial
{\Bbb D}$ then remove one of the 0-handles. This opertaion also removes
an edge of $\partial_f K^\prime_V$. In the following, we generalize this
procedure to all of $K^\prime_I$. Let $L$ be an embedded loop in
$\partial_f K^\prime_I$. We note that $L$ has a natural subdivision into
edges as follows: if $H$ is a 1-handle or 0-handle of $K^\prime_I$ and
$H\cap L$ is non-empty, then $\partial H\cap L = H\cap L$, so that
$H\cap L$ either a single arc or a pair of disjoint arcs. The collection
of these arcs taken over all 0-handles and 1-handles that meet $L$ thus
gives an edge subdivision of $L$. We call an edge of $L$ that is
contained in the boundary of a 0-handle a {\bf 0-handle edge}.
\medskip

\proclaim {Proposition 4.13} Let $f:K\longrightarrow M$ be a tiangular
complex. Assume that $diam(M) > C_3\ell(P)$. Let
$f^\prime:K^\prime\longrightarrow M$ be the complex given by Lemma 4.10. 
Let $K^\prime_I$ denote the complement of the interiors of the $monkey$-handles in
$K^\prime_V$.  Let $L$ be an embedded loop in $\partial_f K^\prime_I$ with edge 
subdivision as above.  If $L$ contains a 0-handle edge, then $L$ contains two 0-handle edges 
and $L$ is inessential in $K^\prime$.
\endproclaim
\medskip

\noindent {\bf Proof} Let ${\Cal A}$ denote the union of all 0-handles 
and 1-handles in $K^\prime_I$ that contain an edge of $L$.  By hypothesis, there is at least
one 0-handle $H$ in ${\Cal A}$ corresponding to the 0-handle edge of $L$.  If $H$  
is joined to another 0-handle by two triangles of $K^\prime$ glued along a common edge, 
then since $L$ is embedded, these are the only two handles in ${\Cal A}$ and 
the proof is comlpete.  So suppose that $H$ is attached to a 1-handle $H_1$ of ${\Cal A}$.  
Homotop $L$ into $H_1$ by pushing the 0-handle edge of $L$ through $H$ into $H_1$.  
Likewise, we can homotop $L$ through $H_1$ into the next handle.  Continuing in this manner,
we must eventually reach a final $0-handle$ disjoint from $H$ since $L$ is embedded.  Contract
$L$ to a point in this 0-handle to complete the proof. $\square$

\medskip
\psfig{file=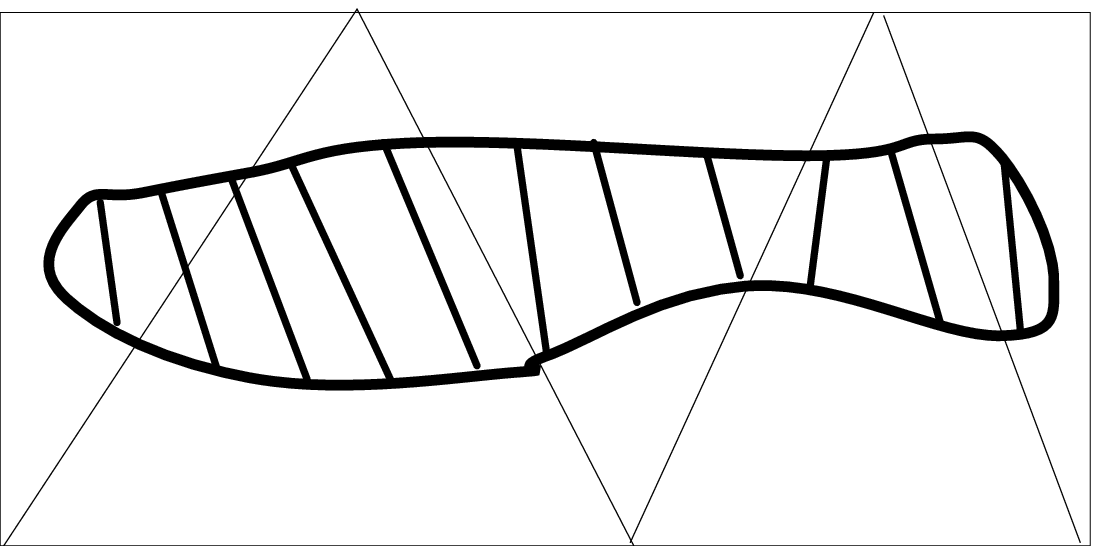,height=2.0in,width=5.0in}
Figure 4.2.
\medskip

\proclaim{Lemma 4.14 (0-handle surgery)} Let $f:K\longrightarrow M$ be a triangular
complex. Assume that $diam(M) > C_3\ell(P)$. Let
$f^\prime:K^\prime\longrightarrow M$ be the complex given by Lemma 4.10. 
Let $K^\prime_I$ denote the I-bundle portion of $K^\prime_V$.  Let $L$ be an embedded loop 
in $\partial_f K^\prime_I$ with edge subdivision as above.  Suppose that $L$ contains a 0-handle
edge and let $H$ denote the corresponding 0-handle. Define $H^+$ to
be the union of $H$ and a small neighborhood of $\partial_f K_I \cap H$ in
the triangle of $K$ containing $H$. There is a 2-complex $K^{\prime\prime}$ with
$K^\prime-int(H^+)\subset K^{\prime\prime}$ and a map 
$f^{\prime\prime}:K^{\prime\prime}-int(H^+)\longrightarrow
M$ is a $\pi_1$ isomorphism. Moreover, $K^{\prime\prime}_V = K^\prime_V - H$.
\endproclaim

\medskip
\noindent {\bf Proof} Attach a disc to $L$ and remove the interior of $H^+$. The result follows
from Lemma 4.1 and Corollary 4.4. $\square$
\medskip
\noindent By applying this result repeatedly, we have the following theorem.

\proclaim {Theorem 4.15 (Good Complex)} Let $f:K\longrightarrow M$ be a triangular
complex. Assume that $diam(M)> C_3\ell(P)$. Let
$f^\prime:K^\prime\longrightarrow M$ be the complex given by Lemma
4.10. There is a $\pi_1$-isomorphic map of a 2-complex 
$f^{\prime\prime}:K^{\prime\prime}\longrightarrow M$ constructed by surgery on 
$K^\prime$ such that, if $K^{\prime\prime}_I$ denotes the I-bundle part of 
$K^\prime_V$, then no simple loop in $\partial_f K^{\prime\prime}_I$
contains a 0-handle edge.  Moreover, 
$Rank(H_1(K^{\prime\prime}_V)) \le B_1\ell(P)^2$ and 
$length(\partial_f K^{\prime\prime}_V) \le 2 length(\partial_f K^\prime_V)$.
\endproclaim
\medskip

\noindent Henceforth, we shall refer to $f^{\prime\prime}:K^{\prime\prime}\longrightarrow M$ as 
constructed above as a {\bf good complex}.  Using this definition, Theorem 4.15 (Good Complex)
gives:
\medskip
\proclaim {Corollary 4.16} Let $f:K\longrightarrow M$ be a good
complex. Assume that $diam(M)> C_3\ell(P)$.  The $\partial_f K_V$ contains at most 
one 0-handle edge for each component of $K_I$.  Hence, there are at most $3\ell(P)$ 
0-handle edges.
\endproclaim
\medskip
\noindent {\bf Proof} If a component of $K_I$ has two attached 0-handles, then there is a loop
in $\partial_f K_I$ containing two 0-handle edges.  This cannot happen in a good complex.  Also, 
there are at most three components of $K_I$ attached to each $monkey$-handle.  Since there are
at most $\ell(P)$ $monkey$-handles, the result follows. $\square$
\medskip

We have now shown that a good complex $K$ has a bounded number of 0-handle edges in 
$\partial_f K_V$.  Since each $monkey$-handle in $K_V$ contributes three edges to 
$\partial_f K_V$, our above discussion shows that we need only bound the number of edges added to
$\partial_f K_V$ by 1-handles.  If fact, no such bound exists; it is possible to have a large 
number of 1-handles which attach very short edges to $K_V$.  The difficulty is that 
$\partial_f K_V$ need not be a closed graph.  Hence it may contain many univalent vertices.  
Instead, we show that there is a closed subgraph, denoted $G_K$, of $\partial_f K_V$ with the 
properties we need.

\medskip
\proclaim {Lemma 4.17} Let $G$ be a closed graph with $Rank(H_1(G))=R$, where $R>1$.  There
is a subdivsion of $G$ into at most $3(R-1)$ edges.
\endproclaim
\medskip
\noindent {\bf Proof} Let $G$ have an arbitrary subdivision into edges and vertices.  Since $G$ 
is closed, there are no univalent vertices.  Given a pair of edges that meet at a bivalent 
vertex, amalgamate these two edges into a single edge, thus removing the bivalent vertex.  
Continue in this manner until every vertex in $G$ has valence at least three.  It then follows
that three times the number of vertices is at most twice the number of edges.  But 
$$R = 1 + (\text{number of edges}) - (\text{number of vertices})$$
whence the result follows. $\square$

\medskip
\psfig{file=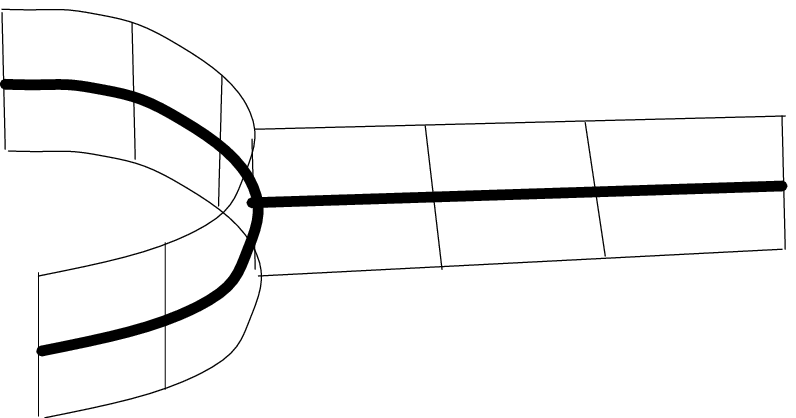,height=2.0in,width=5.0in}
Figure 4.3.The spine $G$ of $K_1$.
\medskip

\proclaim {Lemma 4.18} Let $f:K\longrightarrow M$ be a good complex. Assume that
$diam(M)> C_3\ell(P)$. Then $Rank(H_1(\partial_f K_V))\le B_2\ell(P)^2$ where $B_2 = B_1 + 6$.
\endproclaim
\medskip
\noindent {\bf Proof} Let $K_1$ denote the subcomplex of $K_I$ consisting only of 1-handles.
It follows that $H_1(K_1))=H_1(K_I))$.  There is a spine $G$ of $K_1$ with a subdivision into 
edges and vertices as follows: each 1-handle joins a pair of edges of a triangle in $K$.  
Join these edges with a properly embbeded arc in the 1-handle.  The arc then gives an edge of 
$G$ and its endpoints two vertices of $G$. Continue in this manner while choosing the arcs so 
that their endpoints agree on adjacent 1-handles.  This constructs and subdivides $G$.  
Let $\partial_f K_1 = f^{-1}(f(K_1)\cap \partial V)$.  Now notice that this construction gives a 
natural subdivision of $\partial_f K_1$ into edges an vertices such that $\partial_f K_1$ has 
exactly twice as many edges and exactly twice as many vertices as $G$.  It follows at once
that $H_1(\partial_f K_1))=H_1(G)$ so that in fact $H_1(\partial_f K_1))=H_1(K_I))$.  Now $\partial_f K_V$ is constructed by attaching 0-handle
edges and edges from $monkey$-handles to $\partial_f K_1$.  Each $monkey$-handle attaches three 
edges and Corollary 4.16 shows there are at most $3\ell(P)$ 0-handle edges.  Thus, since there 
are at most $\ell(P)$ $monkey$-handles, we attach at most $6\ell(P)$ edges to $\partial_f K_1$ to
obtain $\partial_f K_V$.  This means that 
$$Rank(H_1(\partial_f K_V)) \le Rank(H_1(G)) + 6\ell(P).$$

\noindent Now using Theorem 4.12 (Bounded Homology) and that $H_1(K_1))=H_1(K_I))$ we have
$$Rank(H_1(\partial_f K_V)) \le \ell(P)^2 + 3\ell(P)$$
\noindent which proves the lemma. $\square$
\medskip

\noindent We comibine Lemma 4.18 and Lemma 4.17 to show that $\partial_f K_V$ has a closed 
subgraph with a bounded number of edges with first homology isomorphic to $H_1(\partial_f K_V)$.
\medskip

\proclaim {Theorem 4.19 (Good Subgraph)} The graph $\partial_f K_V$ contains a closed subgraph $G_K$ with the 
these properties: \newline
\noindent (i) $H_1(G_K) = H_1(\partial_f K_V)$. \newline
\noindent (ii) $G_K$ has a subdivision into at most 
$3B_2\ell(P)^2$ edges. \newline
\noindent  (iii) $length(G_K)\le length(\partial_f K_V)$
\endproclaim

{\bf Proof} Choose a maximal closed subgraph $G_K$ of $\partial_f K_V$.  (i) and (iii) are
obvious and (ii) follows from application of Lemma 4.17 and Lemma 4.18. $\square$
\medskip

\noindent We call a graph $G_K$ provided by Theorem 4.19 a {\bf good subgraph}.

\head 5. Proof Of Main Theorem 
\endhead
\medskip

Suppose the $f:K\longrightarrow M$ is a good complex and assume that 
$diam(M)> C_3\ell(P)$, where $P$ is some presentation of $\pi_1(M)$.
We can prove some results about the homology of loops in 
$\partial_f K_V$.  
To be clear, we have an upper bound on the length of $\partial_f K_V$ 
which precludes (via Corollary 3.8) the existence of embedded loops which 
map to powers of the meridian in
$\partial V$. On the other hand, the Lemma 4.0 shows that there
must exist at least one loop in $\partial_f K_V$ which does map to an
essential loop in $\partial V$. Hence, the image of this essential
loop is also essential in $M$. Let $[\lambda]$ amd $[\mu]$ denote
longitude and meridian basis for $\pi_1(\partial V)$.

\proclaim{Lemma 5.0} Suppose $M$ is a closed, orientable hyperbolic
3-manifold, $P$ is a triangular presentation for $\pi_1(M)$, and
$f:K\longrightarrow M$ is its associated good complex. If
$diam(M)>C_3\ell(P)$ then there exist a pair of embedded loops $\alpha$
and $\beta$ in $\partial_f K_V$ such that both are essential in $K$ and
$\Delta([\alpha], [\beta])\ne 0$. \endproclaim

\noindent {\bf Proof} By Lemma 4.0, there exists at least one such essential loop.
Label this loop $\alpha$. Since $diam(M)>C_3\ell(P)$, we have that
$[\alpha]\ne [\mu]^n$ for any integer $n$. Suppose that
$<[\alpha]>=f_*(H_1(\partial V))$ and let $G$ denote the
component of $K_V$ containing $\alpha$. For simplicity, assume
every other component of $G$ has inessential image in $\partial V$. 
The argument in the proof of Lemma 4.0 shows that we may assume
that each of the other components of $\partial_f K_V$ maps to a
single point in $\partial_f K_V$. We can generate $H_1(G)$ using
$\alpha$ together with a collection of simple closed curves in $G$, none
of which is homologous to $\alpha$ in $G$. The image of each simple
closed curve is null homologous in $H_1(\partial V)$. Thus, we
may modify $f$ homotopically so that each of these simple closed curves
is mapped to a single point. A second homotopy of $f$ makes the image
$f(\alpha)$ embedded. If other components of $G$ have essential image in
$\partial V$, we can repeat the above procedure on each. Since
there are only finitely many components of $G$, we can produce a finite
collection of embedded loops parallel to $f(\alpha)$ which form the
image of $G$. We can then select another loop in $\partial V$
parallel to $f(\alpha)$ that misses $f(K)$, contradicting Lemma 2.3.
$\square$ 

\proclaim{Proposition 5.1} Suppose $A$, $B$, and $\mu$ are closed geodesics
on a Euclidean torus $T$ with $\Delta (A,B) \ne 0$.  Then the 
following inequality holds:
$${1 \over 2} {\ell(\mu) \over \ell(A) + \ell(B)} \le max\lbrace |\Delta
(\mu,A)|, |\Delta (\mu,B)|\rbrace.$$
\endproclaim

\noindent {\bf Proof} Consider the set of lifts of $A$ and $B$ to the
universal cover of $T$.  These form a family of parallelograms 
in the plane.  Each parallelogram has one pair of sides of length
$\ell(A)\over |\Delta(A,B)|$ and another pair of length $\ell(B)\over 
|\Delta(A,B)|$.  Notice that the diameter $D$ of any one of these
parallelograms satisfies
$$\ell(D) \le {{\ell(A) + \ell(B)} \over |\Delta(A,B)|}.$$
It follows that $\ell(\mu) \le [|\Delta(\mu,A)| + |\Delta(\mu,B)|]
\ell(D)$.  One finds immediately that
$${\ell(\mu)\over\ell(A)+\ell(B)} \le {\ell(\mu)\over\ell(A)+\ell(B)} 
|\Delta (A,B)| \le |\Delta(\mu,A)| + |\Delta(\mu,B)|$$
$$\le 2 max(|\Delta(\mu,A)|,|\Delta(\mu,B)|). \square$$

If we combine these two results with with Lemma 3.10 (Short Boundary), we
find that $\partial_f K_V$ contains an embedded loop which is homologous to a huge power of the
core curve of the deep tube.

\proclaim{Corollary 5.2} Suppose $M$ is a closed, orientable hyperbolic
3-manifold, $P$ is a tiangular presentation for $\pi_1(M)$, and
$f:K\longrightarrow M$ is an associated good complex. Let
$[\gamma]$ denote the homotopy class of the core of the deep tube of $M$. 
Let $C>C_3$. If $diam(M)>C\ell(P)$, then there exists an embedded loop $L$ in 
$\partial_f K_V$ such that in $\pi_1(M)$, $f_*([L]) = [\gamma]^N$
where $N \ge {sinh(C\ell(P))\over 4\ell(P)}$.
\endproclaim

\noindent {\bf Proof} Let $\alpha$ and $\beta$ denote the two loops
provided by Lemma 5.1.  It follows from this lemma that
$\Delta(f(\alpha),f(\beta))\ne 0$.  Also, Lemma 3.10 (Short Boundary) gives
that $\ell(f(\alpha)) + \ell(f(\beta)) \le 4\pi\ell(P).$  By applying
Propostion 5.1 to geodesic representatives of $f(\alpha)$ and
$f(\beta)$, we conclude at once that $$max\lbrace |\Delta
(\mu,f(\alpha))|, |\Delta (\mu,f(\beta))|\rbrace \ge {\ell(\mu) \over
8\pi\ell(P)}.$$ Using Lemma 3.3 (Boundary Torus) and $diam(M) \ge C\ell(P)$ we obtain 
$$max\lbrace |\Delta (\mu,f(\alpha))|, |\Delta (\mu,f(\beta))|\rbrace
\ge {sinh(C\ell(P))\over 4\ell(P)}.$$
Now the algebraic intersection of a loop with the merdian equals the
number of times the loop wraps around the core.  This proves the
statement by taking $\gamma$ to be either $\alpha$ or $\beta$ as appropriate. $\square$

It is helpful at this point to arrange that $K_V$ is a connected set.
This can be acheived as follows.  If the loops $X$ and $\beta$ 
provided by Lemma 5.0 are contained in the same component, we agree to 
think of this component as $K_V$.  If they are in different 
components, let $x$ denote a point of $\partial V$ where
$f(X)$ and $f(\beta)$ intersect.  Choose points $y_1$ and $y_2$
on $X$ and $\beta$ respectively which map to $x$.  Join these
points by an gluing the endpoints of an arc to $y_1$ and $y_2$ so that
the arc maps to $f(X)$.  The resulting complex is then connected but
retains the same first homology as the original complex.  In making this
addition, we at most double the length of $\partial_f K_V$
so that we have $length (\partial_f K_V) \le 8\pi e^2\ell (P)$. 

We now wish to utilize some interesting properties which are forced upon
$K_V$ by the geometry in our situation.  In the following, we
prove that we may assume that the map $f:K_V\longrightarrow 
V$ induces an epimorphism between first homology groups.
This is done by showing that the map $f$ lifts to a cover of 
$V$ of index bounded by $\ell(P)$.

\proclaim{Proposition 5.3} Suppose $V$ is a solid torus with Euclidean boundary $T$.
Suppose also that $Area(T) = A$ and $inj(T) = \epsilon$.  Let $\alpha$ and
$\beta$ denote a pair of loops on $T$ such that
$\Delta([\alpha],[\beta])\ne 0$ and $\ell(\alpha)+\ell(\beta) < R$.  
There is a covering $\tilde V$ of index at most ${8R^2\over 3\epsilon^2}$ such that the group 
$<[\tilde \alpha],[\tilde \beta]>$ generates $\pi_1(\tilde V)$.
\endproclaim

\noindent {\bf Proof} Fix a short basis $\lbrace
[X],[Y]\rbrace$ for $\pi_1(T)$.  By Proposition 3.7, we may write
$[\alpha] = a[X]+b[Y]$ and $[\beta] = c[X] + d[Y]$ with
$max\lbrace a,b,c,d\rbrace \le {2R\over \sqrt{3}\epsilon}$.
Let $[\gamma]$ denote the core class of $V$ on $T$.
It follows that $(ad-bc)[\gamma]$ is a finite linear combination of
$[\alpha]$ and $[\beta]$.  This means is $<[\tilde \alpha],[\tilde \beta]>$
isomorphic to $n{\Bbb Z}$ where $n|(ad-bc)$. Thus, by taking an 
at most $|(ad-bc)|$-fold cover of $V$ we obtain the desired lift.  
Moreover, $|(ad-bc)| < 2({4R^2\over 3\epsilon^2})$. $\square$
\medskip

\proclaim{Corollary 5.4} Suppose $M$ is a closed, orientable hyperbolic
3-manifold, $P$ is a triangular presentation for $\pi_1(M)$, and
$f:K\longrightarrow M$ is its associated good complex.
If $diam(M) > C_3\ell(P)$, there is a
covering $\rho:{\Cal Y}\longrightarrow V$ of index at most 
$32(\pi\ell(P))^2\over 3\tilde\epsilon^2$ and a 
lift $\tilde f: K_V\longrightarrow {\Cal Y}$ such that 
$\tilde f_*:\pi_1(K_V)\longrightarrow \pi_1({\Cal Y})$
is an epimorphism.
\endproclaim

\noindent {\bf Proof} Apply Proposition 4.3 using the images under $f$ of 
the loops $\alpha$ and $\beta$ in $\partial_f K_V$ as provided by Lemma 5.1. 
$\square$
\medskip

If we assume that $diam(M) > C_3\ell(P)$, 
we have in effect constructed a 2-complex $f:K_V\longrightarrow V$ which has
very strange properites. First, Theorem 4.12 tells us that the rank of 
$H_1(K_V)$ is bounded by a constant times presentation length.  
Second, we know that the induced homomorphism on first homology is surjective.
Third, by Corollary 5.4, we know that $\partial_f K_V$ contains a short loop 
which maps to a huge power of the core of $V$.  We will now modify $K_V$
to show that this leads to a contradiction.  The first step is to show
that we can attach a small number of discs to $K_V$ to obtain a
new complex which almost gives an $H_1$-isomorphism. This follows by using a good
subgraph of $\partial_f K_V$.
\medskip
\proclaim{Corollary 5.5} Suppose $M$ is a closed, orientable hyperbolic
3-manifold, $P$ is a triangular presentation for $\pi_1(M)$, and
$f:K\longrightarrow M$ is its associated good complex. Let $G_K$ be any good subgraph of 
$\partial_f K_V$. If $diam(M) > C_3\ell(P)$, then the map $\tilde f:G_K\longrightarrow {\Cal Y}$
induces an epimorphism of first homology groups.
\endproclaim
\medskip

\noindent Consider a good subgraph $G_K$.  
Notice that there is a collection of simple closed curves $d_1,...,d_n$ 
such that $\lbrace [d_1],...,[d_n]\rbrace$ is a basis
for $H_1(G_K)$. Now by Lemma 4.19 (Good Subgraph), the length of each loop $d_i$ is 
less than $length(G_K)\le 4\pi (\ell(P))$. The Lemma 4.19 also tells us there are at 
most $B_2\ell(P)^2$ such 
loops since this is the upper bound on the first homolgy rank of 
$H_1(G_K)$. Therefore, the sum of the lengths of these loops is bounded 
above by the product of these two quantities:
$$\sum_{i=1}^n \ell(d_i) \le B_2\ell(P)^2 
\times 4\pi (\ell(P)).$$ 
\noindent We call such a basis for $H_1(G_K)$ a 
{\bf small basis}. $G_K$ is mapped into a torus with bounded geometry via an 
$H_1$-epimorphism.  We use the idea of a small basis to show
that there are restrictions on the complexity of loops which represent 
generators of the kernel of the induced first homolgy map.  

\proclaim {Lemma 5.6} Suppose that $G$ is finite metric graph with
$Rank(H_1(G))\le B$ and $length(G)\le L$.
\newline \noindent Suppose $W$ is a solid torus with
$T=\partial W$ a Euclidean torus with $Area(T) = A$ and $inj(T) = R$.  
Given a continuous map $g:G\longrightarrow W$ such that $g_*$ is an 
epimorphism of first homology groups, there are loops $s_1,...,s_m$ in $G$, 
where $m = rank(H_1(G))-2$ such that $\lbrace [s_1],...[s_m] \rbrace$ is a 
basis for the kernel of $g_*:H_1(G)\longrightarrow H_1(T)$.
Furthermore, for each $1\le i\le m$, $[s_i] = c_{i1}[d_1] +...+
c_{in}[d_n]$ where $[d_1] ... [d_n]$ is a small basis for $H_1(G)$ and
$|c_{ij}| \le ({8 B^2 L^2\over 3R^2})$.
\endproclaim

\noindent {\bf Proof} We have a homomorphism $g_*:H_1(G)\longrightarrow 
H_1(W)$.  We can write $g_*$ as a matrix with integer coefficients

$$g_* = \pmatrix a_{11} && a_{12} && \dots && a_{1n}\cr 
a_{21} && a_{22} &&\dots && a_{2n}
\endpmatrix$$

\noindent using the small basis $\lbrace [d_1],...,[d_n]\rbrace$.  Since
$\sum_{i=1}^n \ell(d_i) \le BL$ we have at once that $\ell(d_i) \le BL$ for 
all $i=1...n$.  Thus, we can use Proposition 3.7 to give an
upper bound on the $|a_{ij}|$.  That is, using a short basis for
$\pi_1(T)$, we conclude that
$|a_{ij}| \le {2BL\over\sqrt(3)R}$ for all $1\le i,j \le n$.  

\noindent Now, since $g_*$ is an epimorphism, this matrix extends to an onto map from 
${\Bbb R^m}$ to ${\Bbb R^2}$.  Hence,
it is not the case that every $2 \times 2$ submatix of $g_*$ has determinant zero.
This means there is a basis for the kernel of $g_*$ using vectors of the form

$$\pmatrix c_{i1} && c_{i2} && 0 && \dots && c_{ij} && 0 && \dots && 0
\endpmatrix$$

\noindent where the entries $c_{ij}$ correspond to those in the cross product of

$$\pmatrix a_{11} && a_{12} && a_{ij} \endpmatrix$$
\noindent with
$$\pmatrix a_{21} && a_{22} && a_{2j} \endpmatrix$$

\noindent It follows that the coefficients $c_{ij}$ are integers that are bounded 
so that for all $1\le i,j \le n$, 
$$|c_{ij}| \le 2max_{ij}(|a_{ij}|)^2\le 2({2BL\over\sqrt(3)R})^2.$$ $\square$
\medskip

Appyling this result to $G_K$ and ${\Cal Y}$ as above, in our situation
we have:

\proclaim {Corollary 5.7} There are loops $s_1,...,s_m$ in $G_K$ 
such that $\lbrace [s_1],...[s_m] \rbrace$ is a basis the kernel of
$f_*:H_1(G_K)\longrightarrow H_1({\Cal Y})$.
Furthermore, for each $1\le i\le m$, $[s_i] = c_{i1}[d_1] +...+
c_{in}[d_n]$, we have $|c_{ij}|\le B_3 \ell(P)^6$ where
$B_3 =  {512\pi^2 B_2^2 \over 3\tilde\epsilon^2}$ 
\endproclaim
\noindent {\bf Proof} Apply Lemma 5.6 with 
$B=B_2\ell(P)^2$, $L = 8\pi\ell(P)$ and 
$R=\tilde\epsilon$. $\square$ 
\medskip

\noindent Finally, we shall need the following algebraic bound in the proof of the
main theorem below:

\proclaim{Proposition 5.8} Suppose $P$ is a presentation of ${\Bbb Z}_N$.
Then $\ell(P) \ge N^{1\over\sqrt{ln (N)}} + \sqrt{ln (N)} - 1$.
\endproclaim
\noindent {\bf Proof} Fix $N > 0$, let $P$ be a presentation of 
${\Bbb Z}_N$.  Since $P$ is a presentation of an Abelian group, 
we have an integer presentation matrix $A = (a_{ij})$ for $P$.  
We can then define $\ell (P) = \Sigma |a_{ij}|$.  We 
may assume that $A$ is a $k \times k$ matrix where $k \le N$.  To see
this, note that if $A$ is not square, we can use 
column and row operations to produce a block presentation matrix $B$ which 
contains a $q\times q$ presentaion matrix in its upper right corner and 
zeros everywhere else.  If $B$ contains more columns than rows, then at 
least one linear combination of the genetators of $P$ is infinite cylcic,
which is clearly impossible since $P$ presents ${\Bbb Z}_N$.  Likewise, 
if $B$ contains more rows than columns, one or more of the rows are
linear combinations of the others.  Hence, we may choose a maximal
collection of linearly independent rows.  By throwing away the remaining
rows, we find the resulting matrix presents a group ${\Bbb Z_S}$ with
$S > N$. Notice that if $A$ is square and $k > N$, then 
$\Sigma |a_{ij}| > N$, so that the $1\times 1$ presentation $(N)$ has 
shorter length than that given by $A$.  All of this means we may assume 
that $A$ is a nonsingular square matrix.  In fact, this implies that
$det(A) = N$.  Now it is well known that 
$$|det(A)| \le \Pi ||A_i||_1$$
\noindent whre $||A_i||_1$ denotes the $L^1$ norm of the $i$th row of
$A$. An easy lower bound for the norm is: 
$$||A||_1 \ge N^{1\over k} + k-1$$
\noindent for every $k\times k$ matrix that presents ${\Bbb Z}_N$. 
This follows since if $||A_i||_1 < N^{1\over k}$ for all $1 \le i\le
k$, then the above inequality gives that $|det(A)|<N$. Since this not
the case, there is at least one element $a_{ij}$ of $A$ with 
$|a_{ij}| \ge N^{1\over k}$. The bound follows since the $k-1$ remaining
rows much each contain an entry of absolute value at least 1. 
 
Fix a value of $N$.  Consider the function $$h(k) = N^{1\over k} +
k-1.$$
\noindent If we minimize this function on the interval $[1,N]$, we find
the absolute minimum is bounded below by $k=\sqrt{ln (N)}$ so
that $$h(k) \ge N^{1\over\sqrt{ln (N)}} + \sqrt{ln (N)} - 1$$ 
\noindent gives a lower bound for every presentation of ${\Bbb Z}_N$.
$\square$

\medskip

\noindent We can now prove the following:

\proclaim {Theorem 5.9} There is an explicit constant $R>0$ such that if
$M$ is a closed, connected, hyperbolic 3-manifold, and $P$ is any 
presentation of its fundamental group, then $diam(M)< R(\ell(P))$.
\endproclaim

{\bf Proof of Main Theorem} We argue by contradiction.  We show that it
is possible to choose $R>C_3$ sufficiently large so that
$diam(M)>R\ell(P)$ is impossible.
We shall contruct the 2-complex $K_V^+$ by
attaching 2-cells to the loops $s_i = c_{i1} d_1+\dots + c_{in} d_n$ 
with the coeffcients $c_{ij}$ provided by Corollary 5.7. The proof of
Theorem 4.12 (Bounded Homology) and our discussions above show that we may 
perform this constrution so that the resulting complex can be triangulated 
with a bounded number of triagles.  To see this,
using the argument in Theorem 4.12 (Bounded Homology), 
retract $K_V$ onto its spine $G$.  In doing this, every $0$-handle
collapses to a vertex of $G$. This graph has the same
properties as the spine of $K_I$, except that we must add a small
neighborhood of a trivalent vertex for each monkey handle.  
Let $T$ denote a maximal tree of $G$.  Crush this tree to a point.  
This produces a $\pi_1$-isomorphic spine $G^\prime$.
Moreover, it follows that from Theorem 4.12 (Bounded Homology) that 
$G^\prime$ has at most $B_1\ell(P)^2$ edges.  
Let $K_V+$ denote the complex built by gluing
discs around the loops which correspond to 
$s_i = c_{i1} d_1+\dots + c_{in} d_n$ as provided by Corollary 5.7.
Hence we can triangulate by coning each of the
attached 2-cells.  The $i$th 2-cell then contains at most 
$$(c_{i1} +\dots + c_{in})\times B_1\ell(P)^2 \le 
B_1\ell(P)^2 \times B_3\ell(P)^6$$

\noindent triangles.  Since the number of attached 2-cells is less than
the rank of $H_1(K_V)$, we obtain from this result and Theorem 4.12
(Bounded Homology)

that $K_V^+$ can be triangulated with at most 
$$B_1^2 B_3\ell(P)^8$$
\noindent triangles.

There is an obvious continuous map $f:K_V^+\longrightarrow Y$ which
induces an epimorphism of first homology groups.  There are essentially
two possibilities for $H_1(K_V^+)$:

\smallskip

\smallskip
\noindent {\bf Case (1)} $H_1(K_V^+)\cong {\Bbb Z}\oplus 
{\Bbb Z}$.  We prove this cannot occur:

\noindent We may form the space $K_V^{++}$ by attaching enough 
discs to $K_V^+$ to make $\pi_1(K_V^{++}) \cong {\Bbb Z}\oplus 
{\Bbb Z}$.  We now consider our original good complex
$f:K\longrightarrow M$.  We can view $K_V^{++}$ as a subset of
$K$.  To see this, notice that we constructed $K_V^{++}$ by
attaching discs to loops that were inessential in $K$.  Thus, by Lemma
4.2 we may assume $KV^{++}$ is contained in $K$.  Then the restrictions 
of the map $f$ induce a pair of homomorphisms $i_1:\pi_1(\partial_f K_{\Cal
V}^{++})\longrightarrow \pi_1(V)$
and $i_2:\pi_1(\partial_f K_{\Cal
V}^{++})\longrightarrow \pi_1(M-int(V))$.  The first map $i_1$ is
injective by our construction.  On the other hand, the incompressibility
of $\partial V$ in $M-int(V)$ ensures that $i_2$ is also
injective.  This is impossible since Van Kampen's Theorem then gives an
injection of $\pi_1(K_V^{++}) \cong {\Bbb Z}\oplus 
{\Bbb Z}$ into $\pi_1(M)$ which contradicts that $M$ is closed
hyperbolic.

\smallskip
\noindent {\bf Case (2)} $H_1(K_V^+)\cong {\Bbb Z}$ or 
$H_1(K_V^+)\cong {\Bbb Z}\oplus {\Bbb Z_r}$ 
\noindent Suppose $H_1(K_V^+)\cong {\Bbb Z}$.  
We show that these cases lead to a contradicition:

Attach a disc to the loop provided by 
Corollary 5.2.  This creates the space $K_V^++$ which has these two 
properties:
\newline\noindent (i) $H_1(K_V^++) = {\Bbb Z}_N$ where   
$N \ge {sinh(R\ell(P))\over 4\ell(P)}$
\newline\noindent (ii) $K_V^++$ has triangulation by $B_4\ell(P)^6$
triangles, where 
$$B_4 = B_2 + B_3$$
\newline 
This constant follows since the loop provided by Corollary 5.2 is
embedded and there are at most $B_2$ edges in $\partial G_K$. 
The triangulation gives a presentation of ${\Bbb Z_N}$ of length at most
$3B_4\ell(P)^6$.  By substituting the value for $N$ in (i) into the 
lower bound expression from Proposition 5.6,
we note that the minimum size of the triangulation in this case is a
function of exponential growth with respect to $R$ and $\ell(P)$ while
the number of triangles exhibits polynomial growth. It then follows that if 
$R>C_3$ is sufficiently large, the construction of 
this complex gives a presentation for ${\Bbb Z}_N$ which is smaller than the lower
bound provided by Proposition 5.6.  A similar argument works if 
$H_1(K_V^+)\cong {\Bbb Z}\oplus {\Bbb Z_r}.$  $\square$

\head References
\endhead
\noindent [A] Adams, Colin C., {\it The Noncompact Hyperbolic 3-Manifold
of Minimal Volume}, Proc. A.M.S. {\bf 100} (1987) 601-606.
\smallskip
\noindent [BGS] Ballman, W.,  Gromov, M., and Schroeder, V.,
{\it Manifolds of Nonpositive Curvature}, Birkhauser (1981).
\smallskip
\noindent [BP] Benedetti, R. and Pettronio, C., {\it Lectures on
Hyperbolic Geometry}, Springer-Verlag (1992).
\smallskip
\noindent [C] Cooper, Daryl, {\it The Volume of a Closed Hyperbolic
3-Manifold is Bounded by $\pi$ Times the Length of any Presentation of
its Fundamental Group}, Proc. A.M.S.
\smallskip
\noindent [Th] Thurston, W.P., {\it The Geometry and Topology of
3-Manifolds}, Princeton University (1979).
\smallskip
\noindent [W] White, M.E., {\it Injectivity Raduius and Fundamental Group
of Hyperbolic 3-manifolds}, To Appear: Communications in Analysis and Geometry.

\enddocument
\bye